\documentclass[12pt]{article}

\usepackage{amssymb,amsmath} 

\DeclareMathAlphabet{\eufrak}{U}{}{}{} 
\SetMathAlphabet\eufrak{normal}{U}{euf}{m}{n}
\SetMathAlphabet\eufrak{bold}{U}{euf}{b}{n}

\numberwithin{equation}{section}

\newenvironment{Proof}{\removelastskip\par\medskip
\noindent{\em Proof.} \rm}{\penalty-20\null\hfill$\square$\par\medbreak}

\newenvironment{Proofy}{\removelastskip\par\medskip
\noindent{\em Proof} \rm}{\penalty-20\null\hfill$\square$\par\medbreak}

 \allowdisplaybreaks

 \def\real{{\mathord{\mathbb R}}}
 
 \def\inte{{\mathord{\mathbb N}}}
 
 \def\qu{{\mathord{\mathbb Z}}}

 \def\Dom{{\mathrm{{\rm Dom}}}}
 \def\trace{{\mathrm{{\rm trace}}}}

 \def\real{{\mathord{{\rm I\kern-3pt R}}}}        
 
 \def\inte{{\mathord{{\rm I\kern-3pt N}}}}
 \def\sZZ{{\rm Z\kern-.45em{}Z}}
 \def\sQQ{{\kern 0.27em \vrule height1.45ex width0.03em depth0em
           \kern-0.30em \rm Q}}
 \def\qu{{\mathchoice
         {\sQQ}
         {\sQQ}
   {\kern 0.225em \vrule height1.05ex width0.025em depth0em \kern-0.25em \rm Q}
   {\kern 0.180em \vrule height0.78ex width0.020em depth0em \kern-0.20em \rm Q}
         }}
 \def\sGG{{\kern 0.27em \vrule height1.45ex width0.03em depth0em
           \kern-0.30em \rm G}}
 \def\gg{{\mathchoice
         {\sGG}
         {\sGG}
   {\kern 0.225em \vrule height1.05ex width0.025em depth0em \kern-0.25em \rm G}
   {\kern 0.180em \vrule height0.78ex width0.020em depth0em \kern-0.20em \rm G}
         }}

 \newtheorem{prop}{Proposition}[section]
 \newtheorem{lemma}[prop]{Lemma}
 \newtheorem{definition}[prop]{Definition}
 \newtheorem{corollary}[prop]{Corollary}
 \newtheorem{theorem}[prop]{Theorem}

 \def\Dom{{\mathrm{{\rm Dom \! \ }}}}
 
 \def\trace{{\mathrm{{\rm trace}}}}

\def\E{\mathop{\hbox{\rm I\kern-0.20em E}}\nolimits}

\newcommand{\dee}{\mbox{$I  \! \! \! \, D$}}
\newcommand{\lee}{\mbox{$I  \! \! \! \, L$}}

 \newcounter{hyp}
\title{\huge Invariance of Poisson measures under random transformations} 
 
\date{February 10, 2011}

\author{
Nicolas Privault
\\ 
Division of Mathematical Sciences 
\\ 
School of Physical and Mathematical Sciences 
\\ 
Nanyang Technological University 
\\ 
SPMS-MAS-05-43, 21 Nanyang Link 
\\ 
Singapore 637371
}

\begin{document}

\hyphenation{func-tio-nals} 
\hyphenation{Privault} 

\maketitle 

\vspace{-1cm} 

\begin{abstract}
{ 
 We prove that Poisson measures are invariant 
 under (random) intensity preserving transformations 
 whose finite difference gradient satisfies a cyclic 
 vanishing condition. 
 The proof relies on moment identities of independent interest 
 for adapted and anticipating Poisson stochastic integrals, 
 and is inspired by the method applied in 
 \cite{ustrotations} on the Wiener space, 
 although the corresponding algebra is 
 more complex than in the Wiener case. 
 The examples of application include transformations 
 conditioned by random sets such as the convex hull of 
 a Poisson random measure. 
} 
\end{abstract}

\baselineskip0.7cm
 
\section{Introduction} 

\noindent 
 Poisson random measures on metric spaces 
 are known to be quasi-invariant 
 under deterministic transformations 
 satisfying suitable conditions, cf. e.g. 
 \cite{vershik2}, \cite{takahashi}. 
 For Poisson processes on the real line this 
 quasi-invariance property also holds under adapted 
 transformations, cf. e.g. \cite{bremaud}, \cite{jacodshiryaev}. 
 The quasi-invariance of Poisson 
 measures on the real line with respect to anticipative transformations 
 has been studied in \cite{girpri} and 
 in the general case of metric spaces in \cite{smorodina1}. 
 In the Wiener case, random non-adapted 
 transformations of Brownian motion have been considered by 
 several authors using the Malliavin calculus, 
 cf. \cite{uzbk} and references therein. 
\\ 
 
 On the other hand, 
 the {\em invariance} property of the law of stochastic processes 
 has important applications, for example to the construction 
 of identically distributed samples of antithetic random variables 
 that can be used for variance reduction in the Monte Carlo method, 
 cf. e.g. \S~4.5 of \cite{brandimarte}. 
 Invariance results for the Wiener measure 
 under quasi-nilpotent random isometries have 
 been obtained in \cite{ustrotations}, 
 \cite{urcr}, 
 by means of the Malliavin calculus, based on the duality between 
 gradient and divergence operators on the Wiener space. 
 In comparison with invariance results, quasi-invariance 
 in the anticipative case 
 usually requires more smoothness on the considered transformation. 
 Somehow surprisingly, the invariance of Poisson 
 measures under non-adapted transformations does not seem 
 to have been the object of many studies to date. 
\\ 

 The classical invariance theorem for Poisson 
 measures states that given a {\em deterministic} transformation 
 $\tau : X \to Y$ between measure spaces $(X,\sigma )$ and 
 $(Y, \mu)$ sending $\sigma$ to $\mu$, the corresponding transformation 
 on point processes maps the Poisson distribution $\pi_\sigma$ 
 with intensity $\sigma (dx)$ on $X$ to the Poisson distribution 
 $\pi_\mu$ with intensity $\mu (dy)$ on $Y$. 
 As a simple deterministic example in the case 
 of Poisson jumps times $(T_k)_{k \geq 1}$ on the half line 
 $X=Y=\real_+$ with $\sigma ( dx ) = \mu ( dx ) = dx/x$, the 
 homothetic transformation $\tau ( x ) = r x$ leaves $\pi_\sigma$ 
 invariant for all fixed $r>0$. 
 However, the random transformation of the Poisson process 
 jump times according to the mapping $\tau ( x ) = x / T_1$ 
 does not yield a Poisson process since the first jump 
 time of the transformed point process is constantly equal 
 to $1$. 
\\ 
 
 In this paper we obtain sufficient conditions for 
 the invariance of {\em random} transformations 
 $\tau : \Omega^X \times X \to Y$ of 
 Poisson random measures on metric spaces $X$, $Y$. 
 Here the almost sure isometry condition on $\real^d$ 
 assumed in the Gaussian 
 case will be replaced by a pointwise condition on the 
 preservation of intensity measures, and the quasi-nilpotence 
 hypothesis 
 will be replaced by a cyclic condition on the finite difference 
 gradient of the transformation, cf. Relation~\eqref{cyclic} 
 below. 
 In particular, this condition is satisfied by predictable 
 transformations of Poisson measures, 
 as noted in Example~\ref{ex1} of Section~\ref{exs}. 
\\ 
 
 In the case of the Wiener space $W = {\cal C}_0 (\real_+ ; \real^d )$ 
 one considers almost surely defined random isometries 
$$ 
 R (\omega ) : L^2 (\real_+;\real^d) \to L^2(\real_+;\real^d), 
 \qquad 
 \omega \in W, 
$$ 
 given by $R (\omega ) h (t) = U ( \omega , t) h(t)$ where 
 $U(\omega , t) : \real^d \to \real^d$, $t\in \real_+$, 
 is a random process of isometries of $\real^d$. 
 The Gaussian character of the measure transformation induced 
 by $R$ is then given by checking for the Gaussianity of the (anticipative) 
 Wiener-Skorohod integral $\delta ( Rh )$ 
 of $Rh$, for all $h\in L^2(\real_+;\real^d)$. 
 In the Poisson case we consider random isometries 
$$ 
 R (\omega ) : L^2_\mu (Y) \to L^2_\sigma (X) 
$$ 
 given by $R (\omega ) h (x) = h(\tau (\omega , x))$ where 
 $\tau (\omega , \cdot ) : ( X , \sigma ) \to ( Y , \mu )$ 
 is a random transformation that maps $\sigma (dx)$ to $\mu (dy)$ 
 for all $\omega \in \Omega^X$. 
 Here, the Poisson character of the measure transformation induced 
 by $R$ is obtained by showing that the 
 Poisson-Skorohod integral $\delta_\sigma ( Rh )$ of $Rh$ 
 has same distribution under $\pi_\sigma$ as the compensated 
 Poisson stochastic integral $\delta_\mu (h)$ of $h$ under $\pi_\mu$, 
 for all $h \in {\cal C}_c ( Y )$. 
\\ 
 
 For this we will use the Malliavin calculus under Poisson measures, 
 which relies on a finite difference gradient $D$ 
 and a divergence operator $\delta$ that extends the 
 Poisson stochastic integral. 
 Our results and proofs are to some extent 
 inspired by the treatment of the Wiener case 
 in \cite{ustrotations}, see \cite{primoment} 
 for a recent simplified proof on the Wiener space. 
 However, the use of {\em finite difference} operators instead 
 of {\em derivation} operators as in the continuous case 
 makes the proofs and arguments 
 more complex from an algebraic point of view. 
\\ 
 
 As in the Wiener case, we will characterize probability measures 
 via their moments. 
 Recall that the moment $E_\lambda [ Z^n ]$ of order $n$ of a Poisson 
 random variable $Z$ with intensity $\lambda$ can be written as 
$$ 
 E_\lambda [ Z^n ] 
 = 
 T_n ( \lambda ) 
$$ 
 where $T_n (\lambda )$ is the Touchard polynomial of order $n$, 
 defined by $T_0 (\lambda )=1$ and the recurrence relation 
\begin{equation} 
\label{tc} 
 T_{n+1} ( \lambda ) 
 = 
 \lambda 
 \sum_{k=0}^{n} 
 {n \choose k} 
 T_k ( \lambda ) 
, 
 \qquad
 n \geq 0, 
\end{equation} 
 also called the exponential polynomials, 
 cf. e.g. \S11.7 of \cite{charalambides}, 
 Replacing the Touchard polynomial $T_n (\lambda )$ by 
 its centered version $\tilde{T}_n (\lambda )$ defined 
 by $\tilde{T}_0 (\lambda )=1$ and 
\begin{equation} 
\label{tcc} 
 \tilde{T}_{n+1} ( \lambda ) 
 = 
 \lambda 
 \sum_{k=0}^{n-1} 
 {n \choose k} 
 \tilde{T}_k ( \lambda ) 
, 
 \qquad 
 n \geq 0, 
\end{equation} 
 yields the moments of the {\em centered} 
 Poisson random variable with intensity $\lambda >0$ 
 as 
$$ 
 \tilde{T}_n ( \lambda ) 
 = 
 E_\lambda [ ( Z - \lambda )^n] 
, 
 \qquad
 n \geq 0. 
$$ 
 Our characterization of Poisson measures will use recurrence relations 
 similar to \eqref{tcc}, cf. \eqref{deltarec} below, and 
 identities for the moments of compensated Poisson stochastic integrals 
 which are another motivation for this paper, cf. Theorem~\ref{l22} 
 below. 
\\ 
 
 The paper is organized as follows. 
 The main results (Corollary~\ref{c09} and Theorem~\ref{c10}) 
 on the invariance of Poisson measures are 
 stated in Section~\ref{mr} after recalling the 
 definition of the finite difference gradient $D$ 
 and the Skorohod integral operator 
 $\delta$ under Poisson measures in Section~\ref{sp}. 
 Section~\ref{exs} contains examples of transformations 
 satisfying the required conditions which include the 
 classical adapted case and transformations acting inside 
 the convex hull generated by Poisson random measures, given the 
 positions of the extremal vertices. 
 Section~\ref{ml} contains the moment identities for Poisson stochastic 
 integrals of all orders that are used in this paper, 
 cf. Theorem~\ref{l22}. 
 In Section~\ref{rmi} we prove the main results of Section~\ref{mr} 
 based on the lemmas on moment identities established in Section~\ref{ml}. 
 In the appendix Section~\ref{app2} 
 we prove some combinatorial results that are needed in the 
 proofs. 
 Some of the results of this paper have been presented 
 in \cite{priinvcr}. 
\section{Poisson measures and finite difference operators} 
\label{sp} 
 In this section we recall the construction of Poisson 
 measures, finite difference operators and Poisson-Skorohod 
 integrals, cf. e.g. 
 \cite{dnop} and \cite{privaultbk2} Chapter~6 for reviews. 
 We also introduce some other operators that will be needed 
 in the sequel, cf. Definition~\ref{d1} below. 
\\ 
 
 Let $X$ be a $\sigma$-compact metric space with 
 Borel $\sigma$-algebra ${\cal B}(X)$ and a 
 $\sigma$-finite diffuse measure $\sigma$. 
 Let $\Omega^X$ denote the configuration space on $X$, i.e. 
 the space of at most countable and locally finite 
 subsets 
 of $X$, 
 defined as 
$$\Omega^X = \left\{ 
 \omega = 
 ( x_i )_{i=1}^N \subset X, \ x_i\not= x_j \ 
 \forall i\not= j, \ N \in \inte \cup \{ \infty \} \right\} 
. 
$$ 
 Each element $\omega$ of $\Omega^X$ is identified with the 
 Radon point measure 
$$ 
 \omega = \sum_{i=1}^{\omega (X)} \epsilon_{x_i}, 
$$ 
 where $\epsilon_x$ denotes the Dirac measure at $x\in X$ 
 and $\omega (X) \in \inte \cup \{ \infty \}$ 
 denotes the cardinality of $\omega$. 
 The Poisson random measure $ N (\omega , d x )$ is defined by 
\begin{equation} 
\label{pntp} 
 N (\omega , d x ) = 
 \omega ( dx ) 
 = 
 \sum_{k=1}^{\omega (X)} \epsilon_{x_k} ( d x ) 
, 
 \qquad 
 \omega \in \Omega^X. 
\end{equation} 
 
 The Poisson probability measure $\pi_\sigma$ on $X$ 
 can be characterized as the only probability measure on $\Omega^X$ 
 under which for all compact disjoint subsets $A_1,\ldots ,A_n$ of $X$, 
 $n\geq 1$, the mapping 
$$ 
 \omega \mapsto ( \omega (A_1 ) , \ldots , \omega (A_n)) 
$$ 
 is a vector of independent Poisson distributed random variables 
 on $\inte$ with respective intensities $\sigma (A_1), \ldots , \sigma (A_n)$. 
\\ 
 
 The Poisson measure $\pi_\sigma$ is also characterized by its 
 Fourier transform 
$$ 
 \psi_{\sigma} ( f ) 
 = 
 E_\sigma \left[ 
 \exp \left( 
 i\int_X f(x) ( \omega (dx) - \sigma (dx)) 
 \right) 
 \right], 
 \qquad 
 f\in L^2_\sigma (X ), 
$$ 
 where $E_\sigma$ denotes expectation under $\pi_\sigma$, 
 which satisfies 
\begin{equation} 
\label{df} 
 \psi_\sigma ( f ) 
 = 
 \exp \left( 
 \int_X ( e^{i f(x)} - i f(x) -1) 
 \sigma (dx) 
 \right) 
, 
 \qquad 
 f\in L^2_\sigma (X ) 
, 
\end{equation} 
 where the compensated Poisson stochastic integral 
 $\displaystyle \int_X f(x) ( \omega (dx) - \sigma (dx ) )$ 
 is defined by the isometry 
\begin{equation} 
\label{isop} 
 E_\sigma \left[ 
 \left( 
 \int_X f(x) ( \omega (dx) - \sigma (dx ) ) 
 \right)^2 
 \right] 
 = 
 \int_X | f(x) |^2 \sigma (dx ), 
 \qquad 
 f\in L^2_\sigma (X ). 
\end{equation} 
 
 We refer to \cite{kree}, \cite{yito}, \cite{picard}, 
 for the following definition. 
\begin{definition} 
 Let $D$ denote the finite difference gradient 
 defined 
 as 
\begin{equation} 
\label{prr} 
 D_x F(\omega ) = \varepsilon^+_x F (\omega ) - F( \omega ), 
 \qquad 
 \omega \in \Omega^X, 
 \quad 
 x\in X, 
\end{equation} 
 for any random variable $F:\Omega^X \to \real$, 
 where 
$$
 \varepsilon^+_x F(\omega ) = F(\omega \cup \{ x \} ), 
 \qquad 
 \omega \in \Omega^X, 
 \quad 
 x\in X. 
$$ 
\end{definition} 
 The operator $D$ is continuous on the space 
 $\dee_{2,1}$ defined by the norm 
$$ 
 \Vert F \Vert_{2,1}^2 = \Vert F \Vert_{L^2 ( \Omega^X , \pi_\sigma )}^2 
 + \Vert D F \Vert_{L^2( \Omega^X \times X , \pi_\sigma \otimes \sigma )}^2, 
 \qquad 
 F \in \dee_{2,1} 
. 
$$ 
 
 We refer to Corollary~1 of \cite{picard} 
 for the following definition. 
\begin{definition} 
 The Skorohod integral operator $\delta_\sigma$ 
 is defined on any measurable 
 process $u : \Omega^X \times X \to \real$ 
 by the expression 
\begin{equation} 
\label{ui} 
 \delta_\sigma ( u ) 
 = 
 \int_X u_t ( \omega \setminus \{ t \} ) 
 ( \omega ( dt ) - \sigma (dt) ) 
, 
\end{equation} 
 provided 
$ 
 \displaystyle 
 E_\sigma \left[ \int_X | u ( \omega , t ) | \sigma (dt ) \right] < \infty 
$. 
\end{definition} 
 Relation \eqref{ui} between $\delta_\sigma$ and the Poisson 
 stochastic integral will be used to characterize the 
 distribution of the perturbed configuration points. 
 Note that if $D_t u_t = 0$, $t\in X$, and 
 in particular when applying \eqref{ui} 
 to $u \in L^1_\sigma ( X )$ a deterministic function, we have 
\begin{equation} 
\label{ui2} 
 \delta_\sigma ( u ) 
 = 
 \int_X u ( t ) ( \omega ( d t ) - \sigma (d t ) ) 
\end{equation} 
 i.e. $\delta_\sigma ( u )$ 
 with the compensated Poisson-Stieltjes integral of $u$. 
 In addition if $X = \real_+$ and $\sigma (dt) = \lambda_t dt$ 
 we have 
\begin{equation} 
\label{*9} 
 \delta_\sigma ( u ) = \int_0^\infty u_t (dN_t - \lambda_t dt) 
\end{equation} 
 for all square-integrable predictable processes $(u_t)_{t\in \real_+}$, 
 where $N_t = \omega ([0,t])$, $t\in \real_+$, is a Poisson process with 
 intensity $\lambda_t >0$, cf. e.g. the Example page 518 of 
 \cite{picard}. 
\\ 
 
 The next proposition can be obtained from Corollaries~1 and 5 in 
 \cite{picard}. 
\begin{prop} 
\label{pp} 
 The operators $D$ and $\delta_\sigma$ are closable 
 and satisfy the duality relation 
\begin{equation} 
\label{dr1} 
 E_\sigma [ \langle D F , u \rangle_{L^2_\sigma ( X )}] 
 = 
 E_\sigma [ F \delta_\sigma ( u ) ], 
\end{equation} 
 on their $L^2$ domains 
 $\Dom ( \delta_\sigma ) 
 \subset L^2 ( \Omega^X \times X , \pi_\sigma \otimes \sigma)$ 
 and 
 $\Dom ( D ) = \dee_{2,1} 
 \subset L^2 ( \Omega^X , \pi_\sigma )$ 
 under the Poisson measure $\pi_{\sigma}$ with intensity $\sigma$. 
\end{prop} 
 
 The operator $\delta_\sigma$ is continuous on 
 the space $\lee_{2,1} \subset \Dom ( \delta_\sigma )$ 
 defined by the norm 
$$ 
 \Vert u \Vert_{2,1}^2 
 = 
 E_\sigma 
 \left[  
 \int_X | u_t |^2 \sigma (d t ) 
 \right] 
 + 
 E_\sigma 
 \left[  
 \int_X | D_s u_t |^2 \sigma (d s ) \sigma (d t ) 
 \right] 
, 
$$ 
 and for any $u\in \lee_{2,1}$ we have the Skorohod isometry 
\begin{equation} 
\label{skois} 
 E_\sigma \left[ \delta_\sigma ( u )^2 \right] 
 = 
 E_\sigma \left[ \Vert u \Vert^2_{L^2_\sigma (X)} \right] 
 + 
 E_\sigma \left[ 
 \int_X \int_X D_s u_t D_t u_s 
 \sigma ( ds ) \sigma ( dt ) 
 \right] 
. 
\end{equation} 
 cf. Corollary~4 and pages 517-518 of \cite{picard}. 
\\ 
 
 In addition, from \eqref{ui} we have the commutation relation 
\begin{equation} 
\label{commrel} 
 \varepsilon^+_t \delta_\sigma (u) 
 = \delta_\sigma ( \varepsilon^+_t u) + u_t, 
 \qquad 
 t \in X 
, 
\end{equation} 
 provided $D_t u \in \lee_{2,1}$, $t\in X$. 
\\ 
 
 The moments identities for Poisson stochastic integrals proved in 
 this paper rely on the decomposition stated in the following lemma. 
\begin{lemma} 
\label{l5} 
 Let $u \in \lee_{2,1}$ be such that 
 $\delta_\sigma ( u )^n \in \dee_{2,1}$, 
 $D_t u \in \lee_{2,1}$, $\sigma (dt)$-a.e., 
 and 
$$ 
 E_\sigma \left[ 
 \int_X 
 | u_t |^{n-k+1} 
 | 
 \delta_\sigma ( \varepsilon^+_t u ) 
 |^k 
 \sigma ( dt ) 
 \right] 
 < \infty, 
 \quad 
 E_\sigma \left[ 
 | \delta_\sigma ( u ) |^k 
 \int_X 
 | u_t |^{n-k+1} 
 \sigma ( dt ) 
 \right] 
 < \infty 
, 
$$ 
 $0 \leq k \leq n$. 
 Then we have 
\begin{eqnarray*} 
 E_\sigma [ \delta_\sigma ( u )^{n+1} ] 
 & = & 
 \sum_{k=0}^{n-1} 
 {n \choose k} 
 E_\sigma \left[ 
 \delta_\sigma ( u )^k 
 \int_X 
 u_t^{n-k+1} 
 \sigma ( dt ) 
 \right] 
\\ 
 & & 
 + 
 \sum_{k=1}^n 
 {n \choose k} 
 E_\sigma \left[ 
 \int_X 
 u_t^{n-k+1} 
 ( 
 \delta_\sigma ( \varepsilon^+_t u )^k 
 - 
 \delta_\sigma ( u )^k 
 ) 
 \sigma ( dt ) 
 \right] 
, 
\end{eqnarray*} 
 for all $n\geq 1$. 
\end{lemma}
\begin{Proof} 
 We have, applying \eqref{commrel} to $F = \delta_\sigma (u)^n$, 
\begin{eqnarray*} 
 E_\sigma [ \delta_\sigma ( u )^{n+1}] 
 & = & 
 E_\sigma \left[ 
 \int_X 
 u_t 
 D_t \delta_\sigma ( u ) ^n 
 \sigma ( dt ) 
 \right] 
\\ 
 & = & 
 E_\sigma \left[ 
 \int_X 
 u_t 
 ( 
 ( 
 \varepsilon^+_t \delta_\sigma ( u ) 
 )^n 
 - 
 \delta_\sigma ( u )^n 
 ) 
 \sigma ( dt ) 
 \right] 
\\ 
 & = & 
 E_\sigma \left[ 
 \int_X 
 u_t 
 ( 
 ( u_t + \delta_\sigma ( \varepsilon^+_t u ) 
 )^n 
 - 
 \delta_\sigma ( u )^n 
 ) 
 \sigma ( dt ) 
 \right] 
\\ 
 & = & 
 \sum_{k=0}^{n-1} 
 {n \choose k} 
 E_\sigma \left[ 
 \int_X 
 u_t^{n-k+1} 
 \delta_\sigma ( \varepsilon^+_t u )^k 
 \sigma ( dt ) 
 \right] 
\\ 
  & = & 
 \sum_{k=0}^{n-1} 
 {n \choose k} 
 E_\sigma \left[ 
 \int_X 
 u_t^{n-k+1} 
 \delta_\sigma ( u )^k 
 \sigma ( dt ) 
 \right] 
\\ 
 & & 
 + 
 \sum_{k=1}^n 
 {n \choose k} 
 E_\sigma \left[ 
 \int_X 
 u_t^{n-k+1} 
 ( 
 \delta_\sigma ( \varepsilon^+_t u )^k 
 - 
 \delta_\sigma ( u )^k 
 ) 
 \sigma ( dt ) 
 \right] 
. 
\end{eqnarray*} 
\end{Proof} 
 From Relation~\eqref{ui2} and Lemma~\ref{l5} we find that 
 the moments of the compensated Poisson stochastic integral 
 $\displaystyle \int_X f (t) (\omega (dt) - \sigma (dt) )$ 
 of $\displaystyle f \in \bigcap_{p=1}^{N+1} L^p_\sigma (X)$ satisfy  
 the recurrence identity 
\begin{eqnarray} 
\label{jnm} 
\lefteqn{ 
 E_\sigma \left[ 
 \left( 
 \int_X f (t) (\omega (dt) - \sigma (dt) ) \right)^{n+1} 
 \right] 
} 
\\ 
\nonumber 
 & = & 
 \sum_{k=0}^{n-1} 
 {n \choose k} 
 \int_X 
 f^{n-k+1} (t)
 \sigma ( dt ) 
 E_\sigma \left[ 
 \left( 
 \int_X f(t) (\omega (dt) - \sigma (dt) ) 
 \right)^k 
 \right] 
, 
\end{eqnarray} 
 $n=0,\ldots ,N$, 
 which is analog to Relation~\eqref{tcc} 
 for the centered 
 Touchard polynomials and coincides with \eqref{isop} for $n=1$. 
\\ 
 
 The Skorohod isometry \eqref{skois} shows that $\delta_\sigma$ is 
 continuous on $\lee_{2,1}$, and that its moment of order two of 
 $\delta_\sigma ( u )$ satisfies 
$$ 
 E_\sigma [\delta_\sigma ( u )^2 ] 
 = 
 E_\sigma [ 
 \Vert u \Vert^2_{L^2_\sigma (X)}  
 ] 
, 
$$
 provided 
$$ 
 \int_X \int_X D_s  u_t D_t u_s \sigma ( ds ) 
 \sigma ( dt ) = 0 
, 
$$ 
 as in the Wiener case \cite{ustrotations}. 
 This condition is satisfied when 
$$ 
 D_t u_s D_s u_t = 0, 
 \qquad 
 s,t \in X 
, 
$$ 
 i.e. 
 $u$ is adapted in the sense of e.g. \cite{privaultwu},  Definition~4, 
 or predictable when $X=\real_+$. 
\\ 
 
 The computation of moments of higher orders 
 turns out to be more technical, cf. Theorem~\ref{l22} below, 
 and 
 will be used to characterize the Poisson distribution. 
 From \eqref{jnm}, in order for 
 $\delta_\sigma ( u ) \in L_\sigma^{n+1} ( \Omega^X )$ 
 to have the same 
 moments as the compensated Poisson integral of 
 $f \in \displaystyle \bigcap_{p=2}^{n+1} L_\sigma^p ( X)$, 
 it should satisfy the recurrence relation 
\begin{equation} 
\label{deltarec} 
 E_\sigma [ \delta_\sigma ( u )^{n+1}] 
 = 
 \sum_{k=0}^{n-1} 
 {n \choose k} 
 \int_X 
 f^{n-k+1} (t)
 \sigma ( dt ) 
 E_\sigma \left[ 
 \delta_\sigma ( u )^k 
 \right] 
, 
\end{equation} 
 $n\geq 0$, 
 which is an extension of Relation~\eqref{jnm} to 
 the moments of compensated Poisson 
 stochastic integrals, and characterizes their distribution 
 by Carleman's condition \cite{carleman} 
 when $\displaystyle \sup_{p\geq 1} \Vert f \Vert_{L^p_\sigma (Y)} < \infty$. 
\\ 
 
 In order to simplify the presentation of 
 moment identities for the Skorohod integral $\delta_\sigma$ 
 it will be convenient to use the following symbolic notation 
 in the sequel. 
\begin{definition} 
\label{d1} 
 For any measurable process $u : \Omega^X \times X \to \real$, let 
\begin{equation} 
\label{defdlta} 
 \Delta_{s_0} 
 \cdots  
 \Delta_{s_j} 
 \prod_{p=0}^n 
 u_{s_p} 
 = 
 \sum_{ 
 { 
 \Theta_0 \cup \cdots \cup \Theta_n = \{ s_0 , s_1 , \ldots , s_j \} 
 \atop 
 s_0 \notin \Theta_0 , \ldots , s_j \notin \Theta_j 
 } 
 } 
 D_{\Theta_0} 
 u_{s_0} 
 \cdots 
 D_{\Theta_n} 
 u_{s_n} 
, 
\end{equation} 
 $s_0,\ldots ,s_n \in X$, $0 \leq j \leq n$, 
 where $\displaystyle D_\Theta : = \prod_{s_j\in \Theta} D_{s_j}$ when 
 $\Theta \subset \{s_0, s_1, \ldots , s_j \}$. 
\end{definition} 
 Note that the sum in \eqref{defdlta} includes empty sets. 
 For example we have 
$$ 
 \Delta_{s_0} 
 \prod_{p=0}^n 
 u_{s_p} 
 = 
 u_{s_0} 
 \sum_{ 
 { 
 \Theta_1 \cup \cdots \cup \Theta_n = \{ s_0 \} 
 \atop 
 s_0 \notin \Theta_0 , \ldots , s_j \notin \Theta_j 
 } 
 } 
 D_{\Theta_1} 
 u_{s_1} 
 \cdots 
 D_{\Theta_n} 
 u_{s_n} 
 = 
 u_{s_0} 
 D_{s_0} 
 \prod_{p=1}^n 
 u_{s_p} 
, 
$$ 
 and 
$ 
 \Delta_{s_0} u_{s_0} = 0
$. 
 The use of this notation allows us to rewrite the Skorohod isometry 
 \eqref{skois} as 
$$ 
 E_\sigma [\delta_\sigma ( u )^2 ] 
 = 
 E_\sigma \left[ \int_X u_s^2 \sigma ( ds ) \right] 
 + 
 E_\sigma \left[ 
 \int_X \int_X \Delta_s \Delta_t ( u_s u_t ) 
 \sigma ( ds ) \sigma ( dt ) 
 \right] 
, 
$$ 
 since by definition we have 
$$ 
 \Delta_s \Delta_t ( u_s u_t ) 
 = 
 D_s u_t D_t u_s, 
 \qquad
 s, t \in X 
. 
$$ 
 As a consequence of Theorem~\ref{l22} 
 and Relation~\eqref{xt} of Proposition~\ref{p01} below, 
 the third moment of $\delta_\sigma (u)$ is given by 
\begin{eqnarray} 
\label{igb} 
\lefteqn{ 
 E_\sigma \left[ 
 \delta_\sigma ( u )^3 
 \right] 
 = 
 E_\sigma \left[ 
 \int_X 
 u^3_s 
 \sigma ( ds ) 
 \right] 
 + 
 3 
 E_\sigma \left[ 
 \delta ( u ) 
 \int_{X} 
 u_s^2 
 \sigma ( ds ) 
 \right] 
} 
\\ 
\nonumber 
 & & 
 \! \! \! \! \! \! \! \! \! \! \! 
 + 
 3 
 E_\sigma \left[  \int_{X^3} 
 \Delta_{s_1} \Delta_{s_2} 
 ( 
 u_{s_1} 
 u_{s_2}^2 
 ) 
 \sigma ( ds_1 ) \sigma ( ds_2 ) 
 \right] 
 + 
 E_\sigma \left[  \int_{X^3} 
 \Delta_{s_1} \Delta_{s_2} \Delta_{s_3} 
 ( 
 u_{s_1} 
 u_{s_2} 
 u_{s_3} 
 ) 
 \sigma ( ds_1 ) \sigma ( ds_2 ) \sigma ( ds_3 ) 
 \right] 
, 
\end{eqnarray}  
 cf. \eqref{qas2} and \eqref{qas3} below, which reduces to 
$$ 
 E_\sigma \left[ 
 \delta_\sigma ( u )^3 
 \right] 
 = 
 E_\sigma \left[ 
 \int_X 
 u^3_s 
 \sigma ( ds ) 
 \right] 
 + 
 3 
 E_\sigma \left[ 
 \delta ( u ) 
 \int_{X} 
 u_s^2 
 \sigma ( ds ) 
 \right] 
$$ 
 when $u$ satisfies the cyclic conditions 
$$ 
 D_{t_1} u_{t_2} D_{t_2} u_{t_1} = 0, 
 \quad 
 \mbox{and} 
 \quad 
 D_{t_1} u_{t_2} D_{t_2} u_{t_3} D_{t_3} u_{t_1} = 0, 
 \qquad 
 t_1,\ldots ,t_3 \in X, 
$$ 
 of Lemma~\ref{l12} in the appendix, which shows that \eqref{defdlta} vanishes, 
 see also \eqref{xt2} below for moments of higher orders. 
 When $X=\real_+$, \eqref{cyclic2} is satisfied 
 in particular when $u$ is predictable with respect to the standard 
 Poisson process filtration. 
\section{Main results} 
\label{mr} 
 The main results of this paper are stated in this section 
 under the form of Corollary~\ref{c09} and 
 Theorem~\ref{c10}. 
\\ 
 
 Let $(Y,\mu)$ denote another measure space with associated 
 configuration space $\Omega^Y$ and $\sigma$-finite diffuse 
 intensity measure $\mu (dy)$. 
 Given an everywhere defined measurable random mapping 
\begin{equation} 
\label{mapping} 
 \tau : \Omega^X \times X \rightarrow Y
,
\end{equation} 
 indexed by $X$, 
 let $\tau_* (\omega )$, $\omega \in \Omega^X$, 
 denote the image measure of $\omega$ 
 by $\tau$, i.e. 
\begin{equation} 
\label{llk} 
 \tau_* : \Omega^X \to \Omega^Y 
\end{equation} 
 maps 
$$ 
 \omega  = \sum_{i=1}^{\omega (X)} 
 \epsilon_{x_i} 
 \in \Omega^X 
 \qquad 
 \mbox{to} 
 \qquad 
 \tau_* (\omega)  = \sum_{i=1}^{\omega (X)} 
 \epsilon_{\tau ( \omega , x_i)} \in \Omega^Y 
. 
$$ 
 In other terms, the random mapping $\tau_* : \Omega^X \to \Omega^Y$ 
 shifts each configuration point $x \in \omega$ according to 
 $x \mapsto \tau ( \omega , x )$, and 
 in the sequel we will be interested in finding conditions 
 for $\tau_*:\Omega^X \to \Omega^Y$ to map $\pi_{\sigma}$ 
 to $\pi_{\mu}$. 
 This question is well known to have an affirmative answer 
 when the transformation $\tau : X \to Y$ is deterministic 
 and maps $\sigma$ to $\mu$, as can be checked from the 
 L\'evy-Khintchine representation \eqref{df} of the characteristic 
 function of $\pi_\sigma$. 
 In the random case we will use the moment identity of 
 the next Proposition~\ref{th}, 
 which is a direct application of Proposition~\ref{c1} below with $u=Rh$. 
 We apply the convention that 
 $\displaystyle \sum_{i=1}^0 l_i = 0$, so that 
 $
 \left\{ 
 l_0 , l_1 \geq 0 \ : \ \sum_{i=1}^0 l_i = 0 
 \right\}$ 
 is an arbitrary singleton. 
\begin{prop} 
\label{th} 
 Let $N\geq 0$ and let $R ( \omega ) : L^p_\mu ( Y ) \to L^p_\sigma ( X )$, 
 $\omega \in \Omega^X$, 
 be a random isometry for all $p = 1, \ldots ,N+1$. 
 Then for all 
 $\displaystyle h \in \bigcap_{p=1}^{N+1} L^p_\mu ( Y )$ 
 such that $Rh \in \lee_{2,1}$ is bounded and 
\begin{equation} 
\label{bound} 
 E_\sigma 
 \left[ 
 \int_{X^{a+1}} 
 \Big| 
 \Delta_{s_0} 
 \cdots  
 \Delta_{s_a} 
 \left( 
 \prod_{p=0}^a 
 \left( 
 Rh (s_p)
 \right)^{l_p} 
 \right) 
 \Big| 
 \sigma ( ds_0 ) \cdots \sigma ( ds_a ) 
 \right] 
 < \infty, 
\end{equation} 
 $l_0 + \cdots + l_a \leq N+1$, $l_0 ,\ldots ,l_a \geq 1$, $a\geq 0$, 
 we have 
 $\delta_\sigma ( Rh ) \in L^{n+1} ( \Omega^X , \pi_\sigma )$ 
 and 
\begin{eqnarray*} 
\lefteqn{ 
 E_\sigma [ \delta_\sigma ( R h )^{n+1}] 
 = 
 \sum_{k=0}^{n-1} 
 {n \choose k} 
 \int_Y 
 h^{n-k+1} ( y )
 \mu ( dy ) 
 E_\sigma \left[ 
 \delta_\sigma ( R h )^k 
 \right] 
} 
\\ 
 & & 
 + 
 \sum_{a=0}^n 
 \sum_{j=0}^a 
 \sum_{b = a}^n 
 \sum_{\substack{l_0 + \cdots + l_a = n - b 
 \\ l_0 ,\ldots ,l_a \geq 0 
 \\ 
 l_{a+1},\ldots ,l_b = 0
 } 
 } 
 {a \choose j} 
 C_{\eufrak{L}_a,b}^{l_0,n} 
\\ 
 & & 
 \times 
 \left( 
 \prod_{q=j+1}^b 
 \int_Y 
 h^{1+l_q} ( y ) 
 \mu (d y ) 
 \right) 
 E_\sigma \left[ 
 \int_{X^{j+1}} 
 \Delta_{s_0} 
 \cdots 
 \Delta_{s_j} 
 \left( 
 \prod_{p=0}^j 
 \left( 
 R h ( s_p )  
 \right)^{1+l_p} 
 \right) 
 d\sigma^{j+1} ( \eufrak{s}_j ) 
 \right] 
, 
\end{eqnarray*} 
 $n=0,\ldots ,N$, where 
 $d\sigma^{j+1} ( \eufrak{s}_j ) = \sigma ( ds_0 ) \cdots \sigma ( ds_j )$, 
 $\eufrak{L}_a = (l_1,\ldots ,l_a )$, and 
\begin{equation} 
\label{deff} 
 C_{\eufrak{L}_a,a+c}^{l_0,n} 
 = 
 (-1)^c 
 {n \choose l_0}  
 \sum_{0 = r_{c+1} < \cdots < r_0 = a+c + 1} 
 \prod_{q=0}^c 
 \prod_{p=r_{q+1} + 1 - ( c - q )}^{r_q - 1 - ( c - q )} 
 {l_1 + \cdots + l_p + p + q - 1  
 \choose l_1 + \cdots + l_{p-1} + p + q -1} 
. 
\end{equation} 
\end{prop} 
 As a consequence of Proposition~\ref{th}, 
 if in addition 
 $R ( \omega ) : L^p_\mu ( Y ) \to L^p_\sigma ( X )$ 
 satisfies the condition 
\begin{equation} 
\label{lc000} 
 \int_{X^{j+1}} 
 \Delta_{t_0} 
 \cdots 
 \Delta_{t_j} 
 \left( 
 \prod_{p=0}^j 
 \left( 
 R h ( \omega , t_p )  
 \right)^{l_p} 
 \right) 
 \sigma ( dt_0 ) \cdots \sigma ( dt_j ) 
 = 
 0 
, 
\end{equation} 
 $\pi_\sigma ( \omega )$-a.s. for all $l_0+\cdots + l_j \leq N+1$, 
 $l_0 \geq 1, \ldots , l_j \geq 1$, $j=1,\ldots ,N$, 
 then we have 
\begin{equation} 
\label{eq01} 
 E_\sigma [ \delta_\sigma ( R h )^{n+1}] 
 = 
 \sum_{k=0}^{n-1} 
 {n \choose k} 
 \int_Y 
 h^{n-k+1} ( y)
 \mu ( dy ) 
 E_\sigma \left[ 
 \delta_\sigma ( R h )^k 
 \right] 
, 
\end{equation} 
 $n = 0 , \ldots , N$, 
 i.e. the moments of $\delta_\sigma (Rh)$ satisfy 
 the extended recurrence relation \eqref{jnm} of the Touchard type. 
\\ 
 
 Hence Proposition~\ref{th} and Lemma~\ref{l12} yield the next corollary in 
 which the sufficient condition \eqref{cyclic} is a strengthened version 
 of the Wiener space condition $\trace ( DRh )^n=0$ of 
 Theorem~2.1 in \cite{ustrotations}. 
\begin{corollary} 
\label{c09} 
 Let $R : L^p_\mu ( Y ) \to L^p_\sigma ( X )$ 
 be a random isometry for all $p \in [1,\infty ]$. 
 Assume that $\displaystyle h \in \bigcap_{p=1}^\infty L^p_\mu ( Y )$ 
 is such that 
 $\displaystyle \sup_{p\geq 1} \Vert h \Vert_{L^p_\mu (Y)} < \infty$, 
 and that $Rh$ satisfies \eqref{bound} and the cyclic condition 
\begin{equation} 
\label{cyclic} 
 D_{t_1} R  h (t_2) \cdots D_{t_k} R h (t_1) = 0, 
 \qquad 
 t_1,\ldots ,t_k \in X 
, 
\end{equation} 
 $\pi_\sigma \otimes \sigma^{\otimes k}$-a.e. 
 for all 
 $k \geq 2$. 
 Then, under $\pi_\sigma$, 
 $\delta_\sigma (Rh)$ 
 has same distribution as 
 the compensated Poisson integral $\delta_\mu ( h )$ 
 of $h$ under $\pi_\mu$. 
\end{corollary} 
\begin{Proof} 
 Lemma~\ref{l12} below shows that Condition~\eqref{lc000} 
 holds under \eqref{cyclic} since 
$$ 
 D_s ( R  h (t) )^l = 
 \varepsilon^+_s ( R  h (t) )^l 
 - 
 ( R  h (t) )^l 
 = 
 \sum_{k=1}^l 
 {l \choose k} 
 ( R  h (t) )^{l-k} 
 ( D_s ( R  h (t) ) )^k 
= 
0
, 
$$ 
 $s, t \in X$, $l \geq 1$, 
 hence by Proposition~\ref{th}, 
 Relation~\eqref{eq01} holds for all $n\geq 1$, 
 and this shows by induction from \eqref{deltarec} that under 
 $\pi_\sigma$, $\delta_\sigma ( R h )$ has 
 same moments as $\delta_\mu (h)$ under $\pi_\mu$. 
 In addition, since 
 $\displaystyle \sup_{p\geq 1} \Vert h \Vert_{L^p_\mu ( Y )} < \infty$, 
 Relation~\eqref{eq01} also shows by induction that 
 the moments of $\delta_\sigma (Rh)$ satisfy 
 the bound $E_\sigma [ | \delta_\sigma ( Rh ) |^n ] \leq ( C n )^n$ 
 for some $C>0$ and all $n\geq 1$, hence they 
 characterize its distribution by the Carleman condition 
$$ 
 \sum_{k=1}^\infty ( 
 E_\sigma [ \delta_\sigma ( Rh )^{2n} ] )^{-1/(2n)} 
 = + \infty 
, 
$$ 
 cf. \cite{carleman} and page 59 of \cite{shohat}. 
\end{Proof} 
 We will apply Corollary~\ref{c09} to the random isometry 
 $R : L^p_\mu ( Y ) \to L^p_\sigma ( X )$ is given as 
$$ 
 Rh = h \circ \tau, \qquad h \in L^p_\mu ( Y ), 
$$ 
 where $\tau : \Omega^X \times X \to Y$ is the random transformation 
 \eqref{mapping} 
 of configuration points considered at the beginning of this section. 
 As a consequence we obtain the following 
 invariance result for Poisson measures when $(X,\sigma ) = (Y,\mu )$. 
\begin{theorem} 
\label{c10} 
 Let $\tau : \Omega^X \times X \to Y$ 
 be a random transformation such that $\tau ( \omega , \cdot ) : X \to Y$ 
 maps $\sigma$ to $\mu$ for all $\omega \in \Omega^X$, 
 i.e. 
$$ 
 \tau_* ( \omega , \cdot ) \sigma = \mu, 
 \qquad 
 \omega \in \Omega^X, 
$$ 
 and satisfying the cyclic condition 
\begin{equation} 
\label{cyclic3} 
 D_{t_1} \tau ( \omega , t_2 ) \cdots D_{t_k} \tau ( \omega , t_1 ) = 0, 
 \qquad 
 \forall 
 \omega \in \Omega^X, 
 \quad 
 \forall  
 t_1,\ldots ,t_k \in X, 
\end{equation} 
 for all $k \geq 1$. 
 Then $\tau_* : \Omega^X \to \Omega^Y$  
 maps $\pi_\sigma$ to $\pi_\mu$, 
 i.e. 
$$ 
 \tau_* \pi_\sigma = \pi_\mu
$$ 
 is the Poisson measure with intensity $\mu (dy)$ on $Y$. 
\end{theorem} 
\begin{Proof} 
 We first show that, under $\pi_\sigma$, $\delta_\sigma ( h \circ \tau )$ 
 has same distribution as 
 the compensated Poisson integral $\delta_\mu ( h )$ 
 of $h$ under $\pi_\mu$, for all $h \in {\cal C}_c (Y)$. 
\\ 
 
 Let $(K_r)_{r\geq 1}$ denote an increasing family of compact subsets of $X$ 
 such that $X = \displaystyle \bigcup_{r\geq 1} K_r$, and let 
 $\tau_r : \Omega^X \times X \to Y$ be defined for $r \geq 1$ 
 by 
$$ 
 \tau_r ( \omega , x ) = \tau ( \omega \cap K_r , x ), 
 \qquad 
 x \in X, \ \omega \in \Omega^X 
. 
$$ 
 Letting $R_r h = h \circ \tau_r$ defines a random isometry 
 $R_r : L^p_\mu ( Y ) \to L^p_\sigma ( X )$ for all $p\geq 1$, 
 which satisfies the assumptions of Corollary~\ref{c09}. 
 Indeed we have 
\begin{eqnarray*} 
 D_s R_r h(t) & = & 
 D_s h(\tau_r ( \omega , t ) ) 
\\ 
 & = & 
 {\bf 1}_{K_r}(s) 
 ( 
 h( 
 \tau_r ( \omega , t ) 
 + 
 D_s \tau_r ( \omega , t ) 
 ) 
 - 
 h(\tau_r ( \omega , t ) ) 
 ) 
\\ 
 & = & 
 {\bf 1}_{K_r}(s) 
 ( 
 h( 
 \tau_r ( \omega \cup \{ s \} , t ) 
 ) 
 - 
 h(\tau_r ( \omega , t ) ) 
 ) 
, 
 \qquad 
 s,t \in X, 
\end{eqnarray*} 
 hence \eqref{cyclic3} implies that 
 Condition~\eqref{cyclic} holds, and 
 Corollary~\ref{c09} shows that we have 
\begin{equation} 
\label{tion} 
 E_\sigma \left[ 
 e^{i \lambda 
 \delta_\mu ( h \circ \tau_r ) } 
 \right] 
 = 
 E_\mu \left[ e^{i \lambda \delta_\mu ( h) } \right] 
, 
\end{equation} 
 for all $\lambda \in \real$. 
 Next we note that Condition \eqref{cyclic3} implies that 
\begin{equation} 
\label{xy} 
 D_t \tau_r ( \omega , t ) = 
 0 
, 
 \qquad
 \forall 
 \omega \in \Omega^X, 
 \quad 
 \forall t \in X 
, 
\end{equation} 
 i.e. $\tau_r (\omega ,t)$ does not depend on the presence 
 or absence of a point in $\omega$ at $t$, and in particular, 
$$ 
 \tau_r ( \omega , t ) = 
 \tau_r (\omega \cup \{ t\} , t ) 
, \qquad t \notin \omega 
, 
$$ 
 and 
$$ 
 \tau_r ( \omega , t ) = 
 \tau_r (\omega \setminus \{ t\} , t ) 
, \qquad t \in \omega 
. 
$$ 
 Hence by \eqref{ui2} we have 
\begin{eqnarray*} 
 \delta_\mu ( h ) \circ \tau_{r*} 
 & = & 
 \int_Y h ( y ) ( \tau_{r*} \omega ( dy ) - \mu ( dy ) ) 
\\ 
  & = & 
 \int_X h ( \tau_r ( \omega , x ) ) ( \omega ( d x ) - \sigma ( d x ) ) 
\\ 
  & = & 
 \int_X h ( \tau_r ( \omega \setminus \{ x \} , x ) ) ( \omega ( d x ) - \sigma ( dx ) ) 
\\ 
  & = & 
 \delta_\sigma ( h \circ \tau_r ) 
, 
\end{eqnarray*} 
 and by \eqref{tion} we get 
\begin{eqnarray*} 
\lefteqn{ 
 E_\sigma \left[ 
 \exp 
 \left(  
 i 
 \int_Y h ( y ) ( \tau_{r*} \omega ( dy ) - \mu ( dy ) ) 
 \right) 
 \right] 
} 
\\ 
 & = & 
 E_\sigma \left[ 
 \exp 
 \left(  
 i 
 \int_Y h ( y ) ( \omega ( dy ) - \mu ( dy ) ) 
 \right) 
 \circ \tau_{r*} 
 \right] 
\\ 
 & = & 
 E_\sigma \left[ 
 e^{i 
 \delta_\mu ( h ) \circ \tau_{r*} } 
 \right] 
\\ 
 & = & 
 E_\mu \left[ e^{i  \delta_\mu ( h) } \right] 
\\ 
 & = & 
 E_\mu \left[ 
 \exp \left( 
 i 
 \int_Y h ( y ) ( \omega ( dy ) - \mu ( dy ) ) 
 \right) 
 \right] 
. 
\end{eqnarray*} 
 Next, letting $r$ go to infinity we get 
$$ 
 E_\sigma \left[ 
 \exp 
 \left(  
 i 
 \int_Y h ( y ) ( \tau_{r*} \omega ( dy ) - \mu ( dy ) ) 
 \right) 
 \right] 
 = 
 E_\mu \left[ 
 \exp \left( 
 i 
 \int_Y h ( y ) ( \omega ( dy ) - \mu ( dy ) ) 
 \right) 
 \right] 
$$ 
 for all $h\in {\cal C}_c ( Y )$, hence the conclusion. 
\end{Proof} 
 In Theorem~\ref{c10} above 
 the identity \eqref{cyclic3} is interpreted 
 for $k\geq 2$ by stating that 
 $\omega \in \Omega^X$, and $t_1,\ldots ,t_k \in X$, 
 the $k$-tuples 
$$ 
 ( 
 \tau ( \omega \cup \{ t_1\} , t_2 ), 
 \tau ( \omega \cup \{ t_2\} , t_3 ), 
 \ldots 
 , 
 \tau ( \omega \cup \{ t_{k-1} \} , t_k ) 
 , 
 \tau ( \omega \cup \{ t_k \} , t_1 ) 
 ) 
$$ 
 and 
$$ 
 ( 
 \tau ( \omega , t_2 ), 
 \tau ( \omega , t_3 ), 
 \ldots 
 , 
 \tau ( \omega , t_k ) 
 , 
 \tau ( \omega , t_1 ) 
 ) 
$$ 
 coincide on at least one component 
 n$^o$ $i\in \{ 1,\ldots ,k\}$ in $Y^k$, 
 i.e. 
$ 
 D_{t_i} \tau ( \omega , t_{i+1 \mod k} ) = 0$. 
\section{Examples} 
\label{exs} 
 In this section we consider some 
 examples of transformations satisfying the hypotheses 
 of Section~\ref{mr}, in case $X=Y$ for $\sigma$-finite 
 measures $\sigma$ and $\mu$. 
 Using various binary relations on $X$ we consider successively 
 the adapted case, and transformations 
 that are conditioned by a random set such as the convex hull of 
 a Poisson random measure. 
 Such results are consistent with the fact that given the position 
 of its extremal vertices, a Poisson random measure remains Poisson 
 within its convex hull, cf. the unpublished manuscript \cite{davydov}, 
 see also \cite{zuyev2} for a related use of stopping sets. 
\\ 
 
\refstepcounter{hyp} 
\label{ex1} 
 {\em \ref{ex1}.} 
 First, we remark that if $X$ is endowed with a total binary relation 
 $\preceq$ and if $\tau : \Omega^X \times X \to Y$ is (backward) 
 predictable in the sense that 
\begin{equation} 
\label{fjdks.222} 
 x \preceq y 
 \quad 
 \Longrightarrow 
 \quad 
 D_x \tau ( \omega , y ) = 0, 
\end{equation} 
 i.e. 
\begin{equation} 
\label{erf1.111} 
 \tau ( \omega \cup \{ x \} , y ) 
 = 
 \tau ( \omega , y ), 
 \qquad 
 x \preceq y 
, 
\end{equation} 
 then the cyclic Condition~\eqref{cyclic3} is satisfied, 
 i.e. we have 
\begin{equation} 
\label{cyclic1.111} 
 D_{x_1} \tau ( \omega , x_2 ) 
 \cdots D_{x_k} \tau ( \omega , x_1 ) = 0, 
 \qquad 
 x_1,\ldots ,x_k \in X, 
 \quad 
 \omega \in \Omega^X 
, 
\end{equation} 
 for all $k \geq 1$. 
 Indeed, for all $x_1,\ldots ,x_k \in X$ 
 there exists $i\in \{ 1,\ldots , k\}$ such that 
 $x_i \preceq x_j$, for all $1\leq j \leq k$, 
 hence $D_{x_i} \tau ( \omega , x_j ) = 0$, $1\leq j \leq k$, 
 by the predictability condition \eqref{fjdks.222}, hence 
 \eqref{cyclic1.111} holds. Consequenly, 
 $\tau_* : \Omega^X \to \Omega^Y$ maps $\pi_\sigma$ to $\pi_\mu$ 
 by Theorem~\ref{c10}, provided $\tau ( \omega , \cdot ) : X \to Y$ 
 maps $\sigma$ to $\mu$ for all $\omega \in \Omega^X$. 
\\ 
 
 Such binary relations on $X$ can be defined via an increasing family 
 $(C_\lambda )_{\lambda \in \real}$ of subsets whose reunion is 
 $X$ and such that for all $x\not= y\in X$ there exists 
 $\lambda_x, \lambda_y \in \real$ 
 with 
 $x \in C_{\lambda_x} \setminus C_{\lambda_y}$ and 
 $y \in C_{\lambda_y}$, 
 or 
 $y \in C_{\lambda_y} \setminus C_{\lambda_x}$ and 
 $x \in C_{\lambda_x}$, 
 which is equivalent to 
 $y \preceq x$ or $x \preceq y$, respectively. 
\\ 
 
 This framework includes the classical adaptedness condition 
 when $X$ has the form $X = \real_+ \times Y$. 
 For example, if 
 $X$ and $Y$ are of the form $X = Y = \real_+ \times Z$, 
 consider the filtration $({\cal F}_t)_{t \in \real_+}$, where 
 ${\cal F}_t$ is generated by 
$$ 
 \{ 
 \sigma ( [0,s] \times A ) \ : \ 
 0 \leq s < t, 
 \ A\in {\cal B}_b ( Z) 
 \} 
, 
$$ 
 where ${\cal B}_c ( Z)$ denotes the compact Borel subsets of 
 $Z$. 
 In this case it is well-known that 
 $\omega \mapsto \tau_* \omega$ is Poisson distributed 
 with intensity $\mu$ under $\pi_\sigma$, 
 provided $\tau (\omega , \cdot ) : \real_+ \times Z \to \real_+ \times Z$ 
 is predictable in the sense that 
 $\omega \mapsto \tau ( s , z )$ is ${\cal F}_t$-measurable for 
 all $0\leq s\leq t$, $z \in Z$, cf. e.g. Theorem~3.10.21 of 
 \cite{bichtelerbk}. 
 Here, Condition~\eqref{fjdks.222} holds for the partial order 
\begin{equation} 
\label{rder} 
 ( s , x ) \preceq ( t , y ) \quad \Longleftrightarrow \quad s \geq t, 
\end{equation} 
 on $Z \times \real_+$  
 by taking $C_\lambda = [ \lambda , \infty ) \times X$, $\lambda \in \real_+$, 
 and the cyclic Condition~\eqref{cyclic3} is satisfied 
 when 
 $\tau (\omega , \cdot ) : \real_+ \times Z \to \real_+ \times Z$ 
 is predictable in the sense of \eqref{fjdks.222}. 
\\ 
 
 Next, we consider other examples in which the binary relation 
 $\preceq$ is configuration dependent. 
 This includes in particular transformations of Poisson 
 measures within their convex hull, given the positions of 
 extremal vertices. 
\\ 
 
\refstepcounter{hyp} 
\label{ex2} 
 {\em \ref{ex2}.} 
 Let $X$ have a compact convex closure in $\real^d$, for example 
 $X = \bar{B} (0,1) \setminus \{0 \}$. 
 For all $\omega \in \Omega^X$, 
 let ${\cal C} (\omega )$ denote the convex hull 
 of $\omega$ in $\real^d$ with interior 
 $\dot{\cal C} (\omega )$, and let 
 $\omega_e = \omega \cap 
 ( {\cal C}(\omega ) \setminus \dot{\cal C}(\omega ) )$ 
 denote the extremal vertices of ${\cal C} (\omega )$. 
 Consider a measurable mapping $\tau : \Omega^X \times X \to X$ 
 such that for all $\omega \in \Omega^X$, 
 $\tau ( \omega , \cdot )$ is measure preserving, maps 
 $\dot{{\cal C}}(\omega )$ to $\dot{{\cal C}}(\omega )$, 
 and for all $\omega \in \Omega^X$, 
\begin{equation} 
\label{by} 
 \tau ( \omega , x ) = 
 \left\{ 
 \begin{array}{ll} 
 \tau ( \omega_e , x ), 
 & 
 x\in \dot{{\cal C}}(\omega ), 
\\ 
\\ 
 x , 
 & 
 x\in X\setminus \dot{{\cal C}}(\omega ), 
\end{array} 
\right. 
\end{equation} 
 i.e. the points of $\dot{{\cal C}}(\omega )$ 
 are shifted by $\tau ( \omega , \cdot )$ depending on the positions 
 of the extremal vertices of the convex hull of $\omega$, 
 which are left invariant by $\tau ( \omega , \cdot )$. 
 The next figure shows an example of a transformation that modifies 
 only the interior of the convex hull generated by the random 
 measure, in which the number of points is taken to be finite for simplicity of 
 illustration. 
\\ 
~ 
\\ 
 
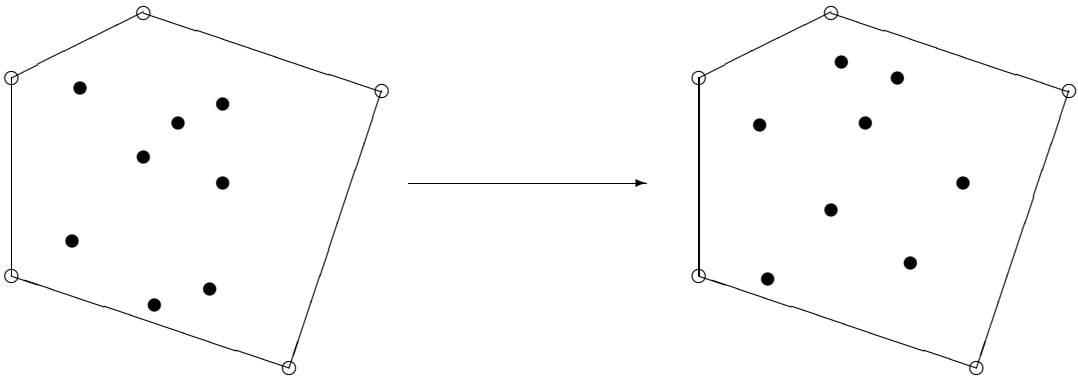
\begin{center} 
\begin{picture}(230,130)(-60,-90)
\linethickness{0.2pt}
\put(-140,40){\circle{5}} 
\put(-140,-35){\circle{5}} 
\put(-140,-35){\line(0,1){75}} 
\put(-140,40){\line(2,1){50}} 
\put(-90,64.5){\circle{5}} 
\put(-90,64.5){\line(3,-1){90}} 
\put(0,35){\circle{5}} 
\put(0,35){\line(-1,-3){35}} 
\put(-35,-70){\circle{5}} 
\put(-35,-70){\line(-3,1){105}} 
\put(-60,0){\circle*{5}} 
\put(-65,-40){\circle*{5}} 
\put(-60,30){\circle*{5}} 
\put(-90,10){\circle*{5}} 
\put(-117,-22){\circle*{5}} 
\put(-77,23){\circle*{5}} 
\put(-86,-46){\circle*{5}} 
\put(-114,36){\circle*{5}} 

\put(10,0){\vector(3,0){90}} 

\put(120,40){\circle{5}} 
\put(120,-35){\circle{5}} 
\put(120,-35){\line(0,1){75}} 
\put(120,40){\line(2,1){50}} 
\put(170,64.5){\circle{5}} 
\put(170,64.5){\line(3,-1){90}} 
\put(260,35){\circle{5}} 
\put(260,35){\line(-1,-3){35}} 
\put(225,-70){\circle{5}} 
\put(225,-70){\line(-3,1){105}} 
\put(220,0){\circle*{5}} 
\put(195,40){\circle*{5}} 
\put(200,-30){\circle*{5}} 
\put(170,-10){\circle*{5}} 
\put(143,22){\circle*{5}} 
\put(183,23){\circle*{5}} 
\put(174,46){\circle*{5}} 
\put(146,-36){\circle*{5}} 
\end{picture}
\end{center} 
\vspace{-0.5cm} 
 
 Next we prove the invariance of such transformations as a consequence 
 of Theorem~\ref{c10}. 
 This invariance property is related to the intuitive fact that 
 given the positions of the extreme vertices, the distribution 
 of the inside points remains Poisson 
 when they are shifted according to the data of the vertices, 
 cf. e.g. \cite{davydov}. 
\\ 
 
 Here we consider the binary relation $\preceq_\omega$ given by 
$$ 
 x \preceq_\omega y 
 \quad 
 \Longleftrightarrow 
 \quad 
 x \in {\cal C} ( \omega \cup \{ y \} ), 
 \qquad 
 \omega \in \Omega^X, 
 \quad x,y \in X 
. 
$$ 
 The relation $\preceq_\omega$ is clearly reflexive, and it is 
 transitive since $x\preceq_\omega y$ and $y\preceq_\omega z$ 
 implies 
$$ 
 x 
 \in 
 {\cal C} ( \omega \cup \{ y \} ) 
 \subset 
 {\cal C} ( \omega \cup \{ z \} ) 
 , 
$$ 
 hence $x\preceq_\omega z$. 
 Note that $\preceq_\omega$ is also total on ${\cal C}(\omega )$ and it 
 is an order relation  on $X \setminus {\cal C} (\omega )$, 
 since it is also antisymmetric on that set, 
 i.e. if $x,y \notin {\cal C} (\omega )$ then 
$$ 
 x\preceq_\omega y \quad 
 {\mbox and} 
 \quad 
 y\preceq_\omega x 
$$ 
 means 
$ 
 x \in {\cal C} ( \omega \cup \{ y \} ) 
$ 
 and 
$ 
 y \in {\cal C} ( \omega \cup \{ x \} ) 
$, 
 which implies $x=y$. 
 We will need the following lemma. 
\begin{lemma} 
\label{ll00} 
 For all $x, y \in X$ and $\omega \in \Omega^X$ 
 we have 
\begin{equation} 
\label{f1} 
 x \preceq_\omega y 
 \quad 
 \Longrightarrow 
 \quad 
 D_x \tau ( \omega , y ) = 0, 
\end{equation} 
 and 
\begin{equation} 
\label{f2} 
 x \not\preceq_\omega y 
 \quad 
 \Longrightarrow 
 \quad 
 D_y \tau ( \omega , x ) = 0. 
\end{equation} 
\end{lemma} 
\begin{Proof} 
 Let $x, y \in X$ and $\omega \in \Omega^X$. 
 First, if $x \not\preceq_\omega y$ then we have 
 $x \notin {\cal C} ( \omega \cup \{ y \} )$ 
 hence 
$ 
 \tau ( \omega \cup \{ y \} , x ) = 
 \tau ( \omega , x ) = x $ by \eqref{by}. 
 Next, if $x \preceq_\omega y$, i.e. $x \in {\cal C} ( \omega \cup \{ y \} )$, 
 we can distinguish two cases: 
\begin{description} 
\item{a)} 
 $x \in {\cal C}(\omega)$. In this case we have 
 ${\cal C}(\omega \cup \{ x \}) = {\cal C}(\omega ) $, 
 hence 
$ 
 \tau ( \omega \cup \{x \} , y ) = \tau ( \omega , y )
$ 
 for all $y \in X$.  
\item{b)} 
 $x \in {\cal C} ( \omega \cup \{ y \} ) \setminus {\cal C} ( \omega )$. 
 If $y \in {\cal C} ( \omega \cup \{ x \} )$ then 
 $x=y \notin \dot{\cal C} ( \omega \cup \{x\} )$, 
 hence 
$ 
 \tau ( \omega \cup \{x \} , y ) = \tau ( \omega , y )
$. 
 On the other hand if 
 $y \notin {\cal C} ( \omega \cup \{ x \} )$ 
 then $y \not\preceq_\omega x$ and 
$ 
 \tau ( \omega \cup \{x \} , y ) = 
 \tau ( \omega , y ) = y
$ as above. 
\end{description} 
 We conclude that $D_x \tau ( \omega , y ) = 0$ 
 in both cases. 
\end{Proof} 
 Let us now show that $\tau : \Omega^X \times X \to \Omega^X$ 
 satisfies the cyclic condition \eqref{cyclic3}. 
 Let $t_1,\ldots ,t_k\in X$. 
 First, if $t_i \in {\cal C} ( \omega )$ 
 for some $i \in \{1,\ldots ,k\}$, then for all $j=1,\ldots ,k$ 
 we have $t_i \preceq_\omega t_j$ and by Lemma~\ref{ll00} we get 
$$ 
 D_{t_i} \tau ( \omega , t_j ) = 0, 
$$ 
 thus \eqref{cyclic3} holds, and 
 we may assume that $t_i \notin {\cal C} ( \omega )$ 
 for all $i = 1,\ldots ,k$. 
 In this case, if $t_{i+1 \! \! \mod k} \not\preceq_\omega t_i$ 
 for some $i=1,\ldots ,k$, then by Lemma~\ref{ll00} we have 
$$ 
 D_{t_i} \tau ( \omega , t_{i+1 \! \! \mod k} ) = 0 
, 
$$ 
 which shows that \eqref{cyclic3} holds. 
 Finally, if 
 $t_1 \preceq_\omega t_k \preceq_\omega \cdots \preceq_\omega t_2 
 \preceq_\omega t_1$, 
 then by transitivity of $\preceq_\omega$ we have 
 $t_1 \preceq_\omega t_k \preceq_\omega t_1$, 
 which implies $t_1=t_k \notin {\cal C} ( \omega )$ by antisymmetry 
 on $X \setminus {\cal C} (\omega )$, 
 hence $D_{t_k} \tau ( \omega , t_1 ) = 0$, 
 and ${\tau} : \Omega^X \times X \to X$ 
 satisfies the cyclic Condition~\eqref{cyclic3} 
 for 
 all $k\geq 2$. 
 Hence $\tau$ satisfies the hypotheses of Theorem~\ref{c10}, 
 and $\tau_* \pi_\sigma = \pi_\mu$ provided 
 $\tau ( \omega , \cdot ) : X \to Y$ 
 maps $\sigma$ to $\mu$ for all $\omega \in \Omega^X$. 
\section{Moment identities for stochastic integrals} 
\label{ml} 
 In this section we prove a moment identity for 
 Poisson stochastic integrals of arbitrary orders in 
 Theorem~\ref{l22}, 
 whose application will be to prove Proposition~\ref{th}. 
 More precisely, given $F : \Omega^X \to \real$ a random variable and 
 $u : \Omega^X \times X \to \real$ a measurable process, 
 we aim at decomposing 
$ 
 E_\sigma \left[ 
 \delta_\sigma ( u )^n 
 F 
 \right] 
$ 
 in terms of the gradient $D$, 
 while removing all occurrences of $\delta_\sigma$ 
 using the integration by parts formula \eqref{dr1}. 
\\ 
 
 In Theorem~\ref{l22} and in the rest of this section 
 we will use the notation 
$$ 
 \varepsilon^+_{\eufrak{s}_b} 
 = 
 \varepsilon^+_{s_1} 
 \cdots 
 \varepsilon^+_{s_b} 
, 
 \qquad 
 \eufrak{s}_b = (s_1,\ldots ,s_b ) \in X^b, 
 \quad 
 b \geq 1. 
$$ 
 Moreover, by saying that 
 $u:\Omega^X \times X \to \real$ has a compact support in $X$ we 
 mean that there exists a compact subset $K$ of $X$ such that 
 $u(\omega , x ) = 0$ for all $\omega \in \Omega^X$ and $x \in X \setminus K$. 
\begin{theorem} 
\label{l22} 
 Let $F:\Omega^X \to \real$ be a bounded random variable and 
 let $u:\Omega^X \times X \to \real$ be a bounded process 
 with compact support in $X$. For all $n\geq 0$ we have 
\begin{equation} 
\label{plkmn} 
 E_\sigma \left[ 
 \delta_\sigma ( u )^n 
 F 
 \right] 
 = 
 \sum_{a=0}^n 
 \sum_{b = a}^n 
 (-1)^{b-a} 
 \! \! \! \! \! \! \! \! \! 
 \sum_{ 
 \substack{ 
 l_1 + \cdots + l_a = n - b 
 \\ 
 l_1,\ldots ,l_a \geq 0 
 \\ 
 l_{a+1},\ldots ,l_b = 0 
 } 
 } 
 \! \! \! \! \! \! \! \! \! 
 C_{\eufrak{L}_a,b} 
 E_\sigma \left[  \int_{X^b} 
 \! \! \! 
 \varepsilon^+_{\eufrak{s}_a } 
 F 
 \prod_{p=1}^b 
 \varepsilon^+_{ \eufrak{s}_a \setminus s_p} 
 u_{s_p}^{1+l_p} 
 \hskip0.05cm 
 d\sigma^b ( \eufrak{s}_b ) 
 \right] 
, 
\end{equation} 
 where 
 $d\sigma^b ( \eufrak{s}_b ) = \sigma ( ds_1 ) \cdots \sigma ( ds_b )$, 
 $\eufrak{L}_a = (l_1,\ldots ,l_a )$, and 
\begin{equation} 
\label{defin} 
 C_{\eufrak{L}_a,a+c} 
 = 
 \! \! \! \! \! \! \! 
 \sum_{0 = r_{c+1} < \cdots < r_0 = a+c+1} 
 \ 
 \prod_{q=0}^c 
 \ 
 \prod_{p=r_{q+1} + q - c + 1 }^{r_q+q-c-1} 
 {l_1 + \cdots + l_p + p + q - 1  
 \choose l_1 + \cdots + l_{p-1} + p + q - 1 } 
. 
\end{equation} 
\end{theorem} 
 Before turning to the proof of Theorem~\ref{l22} we 
 consider some examples. 
\vskip0.2cm  
{\em 1.} 
 For $n=2$ and $F=1$, Theorem~\ref{l22} recovers 
 the Skorohod isometry \eqref{skois} as follows: 
\small 
\begin{align} 
\nonumber 
 E_\sigma \left[ 
 \delta_\sigma ( u )^2 
 \right] 
 = \ & 
 E_\sigma \left[  \int_{X^2} 
 u_{s_1} u_{s_2} 
 \sigma ( ds_1 ) \sigma ( ds_2 ) 
 \right] 
 & [a=0,b=2] 
\\ 
\nonumber 
 & 
 - 
 2 
 E_\sigma \left[  \int_{X^2} 
 u_{s_1} 
 ( I + D_{s_1} ) 
 u_{s_2} 
 \sigma ( ds_1 ) \sigma ( ds_2 ) 
 \right] 
 & [a=1,b=2] 
\\ 
\nonumber 
 & 
 + 
 E_\sigma \left[  
 \int_X 
 | u_{s_1} |^2 
 \sigma ( ds_1 ) 
 \right] 
 & 
 [ a=1,b=1] 
\\ 
\nonumber 
 & 
 + 
 E_\sigma \left[  \int_{X^2} 
 ( I + D_{s_1} ) 
 u_{s_2} 
 ( I + D_{s_2} ) 
 u_{s_1} 
 \sigma ( ds_1 ) \sigma ( ds_2 ) 
 \right] 
 & 
 [ a=2 , b=2] 
\\ 
\label{qas1} 
 = & \ 
 E_\sigma \left[  
 \int_X 
 | u_s |^2 
 \sigma ( ds ) 
 \right] 
 + 
 E_\sigma \left[  \int_{X^2} 
 \Delta_{s_1} 
 \Delta_{s_2} 
 ( 
 u_{s_1} 
 u_{s_2} 
 ) 
 \sigma ( ds_1 ) \sigma ( ds_2 ) 
 \right] 
. 
 & 
\end{align} 
{\em 2.} For $n=3$ and $F=1$, Theorem~\ref{l22} yields the 
 following third moment identity: 
\small 
\begin{align} 
\nonumber 
 & 
 E_\sigma \left[ 
 \delta_\sigma ( u )^3 
 \right] 
 = 
 E_\sigma \left[ 
 \int_X 
 u^3_{s_1} 
 \sigma ( ds_1 ) 
 \right] 
 & 
 [a=1,b=1] 
\\ 
\nonumber 
 & 
 - 
 3 
 E_\sigma \left[  \int_{X^2} 
 u^2_{s_1} 
 ( I + D_{s_1} ) 
 u_{s_2} 
 \sigma ( ds_1 ) \sigma ( ds_2 ) 
 \right] 
 & 
 [a=1,b=2] 
\\ 
\nonumber 
 & 
 + 
 3 
 E_\sigma \left[  \int_{X^2} 
 ( I + D_{s_2} ) 
 u_{s_1} 
 ( I + D_{s_1} ) 
 u_{s_2}^2 
 \sigma ( ds_1 ) \sigma ( ds_2 ) 
 \right] 
 & 
 [ a=2 , b=2] 
\\ 
\nonumber 
 & 
 - 
 E_\sigma \left[  \int_{X^3} 
 u_{s_1} u_{s_2} u_{s_3} 
 \sigma ( ds_1 ) \sigma ( ds_2 ) \sigma ( ds_3 ) 
 \right] 
 & 
 [a=0,b=3] 
\\ 
\nonumber 
 & 
 + 
 3 
 E_\sigma \left[  \int_{X^3} 
 u_{s_1} 
 ( I + D_{s_1} ) 
 u_{s_3} 
 ( I + D_{s_1} ) 
 u_{s_2} 
 \sigma ( ds_1 ) \sigma ( ds_2 ) \sigma ( ds_3 ) 
 \right] 
 & 
 [a=1,b=3] 
\\ 
\nonumber 
 & 
 - 
 3 
 E_\sigma \left[  
 \int_{X^3} 
 (I+D_{s_1}) (I+D_{s_2}) 
 u_{s_3} 
 (I+D_{s_1}) 
 u_{s_2} 
 (I+D_{s_2}) 
 u_{s_1} 
 \right. 
 & 
\\ 
\nonumber 
 & 
 \ \ \ 
 \sigma ( ds_1 )  \sigma ( ds_2 )  \sigma ( ds_3 ) 
 \Big] 
 & 
 [ a=2,b=3] 
\\ 
\nonumber 
 & 
 + 
 E_\sigma \left[  \int_{X^3} 
 ( I + D_{s_1} ) ( I + D_{s_2} ) 
 u_{s_3} 
 ( I + D_{s_1} ) ( I + D_{s_3} ) 
 u_{s_2} 
 ( I + D_{s_2} ) ( I + D_{s_3} ) 
 u_{s_1} 
\right. 
& 
\\ 
\nonumber 
 & 
 \ \ \ 
 \sigma ( ds_1 ) \sigma ( ds_2 ) \sigma ( ds_3 ) 
 \Big] 
 & 
 [ a=3 , b=3] 
\end{align} 
\vskip-0.6cm 
\begin{align} 
\nonumber 
 = & \ 
 E_\sigma \left[ 
 \int_X 
 u^3_{s_1} 
 \sigma ( ds_1 ) 
 \right] 
 + 
 3 
 E_\sigma \left[  \int_{X^2} 
 u_{s_1} 
 D_{s_1} 
 u_{s_2}^2 
 \sigma ( ds_1 ) \sigma ( ds_2 ) 
 \right] 
 & 
\\ 
\label{qas2} 
 & 
 \! \! \! \! \! \! 
 + 
 3 
 E_\sigma \left[  \int_{X^3} 
 \Delta_{s_1} \Delta_{s_2} 
 ( 
 u_{s_1} 
 u_{s_2}^2 
 ) 
 \sigma ( ds_1 ) \sigma ( ds_2 ) 
 \right] 
 + 
 E_\sigma \left[  \int_{X^3} 
 \Delta_{s_1} \Delta_{s_2} \Delta_{s_3} 
 ( 
 u_{s_1} 
 u_{s_2} 
 u_{s_3} 
 ) 
 \sigma ( ds_1 ) \sigma ( ds_2 ) \sigma ( ds_3 ) 
 \right] 
. 
 & 
\end{align} 
\normalsize 

\baselineskip0.7cm 
{\em 3.} 
 Noting that $C_{\eufrak{L}_a,c}$ defined in \eqref{defin} represents 
 the number of partitions of a set of $l_1+\cdots + l_a+a+c$ elements 
 into $a$ subsets of lengths $1+l_1,\ldots ,1+l_a$ and $c$ singletons, 
 we find that when $F=1$ and $u = {\bf 1}_A$ is a deterministic 
 indicator function, and Theorem~\ref{l22} reads 
$$ 
 E_\sigma \left[ 
 ( Z - \lambda )^n 
 \right] 
 = 
 \sum_{a=0}^n 
 \lambda^a 
 \sum_{c=0}^a 
 (-1)^c 
 {n \choose c} 
 S ( n - c , a - c ) 
$$
 for $Z - \lambda = \delta ( {\bf 1}_A ) = \omega ( A ) - \sigma (A)$ 
 a compensated Poisson random variable with intensity 
 $\lambda = \sigma (A)$, where $S(n,c)$ denotes the Stirling number of the 
 second kind, i.e. the number of ways to partition a set of $n$ objects 
 into $c$ non-empty subsets. 
 This coincides with the moment formula 
$$ 
 E_\lambda \left[ 
 ( Z - \lambda )^n 
 \right] 
 = 
 \sum_{a=0}^n 
 \lambda^a 
 S_2(n,a) 
, 
$$ 
 where $S_2(n,a)$ denotes the number of partitions 
 of a set of size $n$ into $a$ non-singleton subsets, 
 which can be obtained from the sequence 
 $(0,\lambda , \lambda, \ldots )$ 
 of cumulants of the compensated Poisson distribution, 
 through the combinatorial identity 
$$ 
 S_2 ( n , a ) 
 = 
 \sum_{c = 0}^a 
 (-1)^c 
 {n \choose c } 
 S ( n - c , a - c ), \qquad 0 \leq a \leq n, 
$$ 
 which is the binomial dual of 
$$ 
 S ( m , n ) 
 = 
 \sum_{k=0}^n 
 {m \choose k} 
 S_2 ( m - k , n - k ) 
.
$$ 
 The proof of Theorem~\ref{l22} will be done by induction based 
 on the following lemma. 
\begin{lemma} 
\label{l11} 
 Let $G :\Omega^X \to \real$ be a bounded random variable and 
 let $u:\Omega^X \times X \to \real$ be a bounded process 
 with compact support in $X$. 
 For all $n\geq 0$ we have 
\begin{eqnarray} 
\label{1.1.1} 
\lefteqn{ 
 E_\sigma \left[ 
 \delta_\sigma ( u )^n 
 G 
 \right] 
 = 
 \sum_{0 = k_d < \cdots < k_0 = n \atop 0 \leq d \leq n} 
 c_{ \eufrak{K}_d } 
 E_\sigma \left[  \int_{X^d} 
 \varepsilon^+_{\eufrak{s}_d} G 
 \prod_{p=1}^d 
 \varepsilon^+_{\eufrak{s}_d \setminus s_p} u_{s_p}^{k_{p-1} - k_p} 
 \hskip0.05cm 
 d\sigma^d ( \eufrak{s}_d ) 
 \right] 
} 
\\ 
\nonumber 
 & & 
\! \! \! \! \! \! \! \! \! \! 
 - 
 \sum_{0 = k_d < \cdots < k_0 = n \atop 1 \leq d \leq n} 
 \! \! \! \! \! 
 c_{ \eufrak{K}_d } 
 E_\sigma \left[ 
 \int_{X^d} 
 \delta_\sigma (  \varepsilon^+_{\eufrak{s}_{d-1}} 
 u 
 )^{k_{d-1}-1}
 \varepsilon^+_{\eufrak{s}_{d-1}} 
 ( 
 u_{s_d} 
 G 
 ) 
 \prod_{p=1}^{d-1} 
 \varepsilon^+_{\eufrak{s}_{d-1} \setminus s_p} 
 u_{s_p}^{k_{p-1} - k_p} 
 \hskip0.05cm 
 d\sigma^d ( \eufrak{s}_d ) 
 \right] 
, 
\end{eqnarray} 
 where 
$ 
\displaystyle 
 c_{\eufrak{K}_d} 
 = 
 \prod_{p=0}^{d-1} 
 {k_p - 1 \choose k_{p+1}} 
$, 
 $\eufrak{K}_d = (k_0,\ldots ,k_d ) \in \inte^{d+1}$. 
\end{lemma} 
\begin{Proof} 
 The formula clearly holds when $n=0$, while 
 when $n\geq 1$, the first summation in \eqref{1.1.1} 
 actually starts from $d=1$. 
 The proof follows by application to $l=n-1$ or $l=n$ 
 of the following identity: 
\begin{eqnarray} 
\label{1n2} 
\lefteqn{ 
 E_\sigma \left[ 
 \delta_\sigma ( u )^n 
 G 
 \right] 
} 
\\ 
\nonumber 
 & = & 
 \sum_{0 = k_{l+1} < \cdots < k_0 = n} 
 c_{\eufrak{K}_{l+1}} 
\\ 
\nonumber 
 & & 
 \! \! \! \! \! \! \! 
 E_\sigma \left[ 
 \int_{X^{l+1}} 
 \delta_\sigma ( 
 \varepsilon^+_{\eufrak{s}_l } 
 u 
 )^{k_l-1}
 \varepsilon^+_{\eufrak{s}_{l+1}} 
 G 
 \varepsilon^+_{\eufrak{s}_l} u_{s_{l+1}} 
 \varepsilon^+_{s_{l+1} } 
 \prod_{p=1}^l 
 \varepsilon^+_{\eufrak{s}_{l+1} \setminus s_p} 
 u_{s_p}^{k_{p-1} - k_p} 
 \hskip0.05cm 
 d\sigma^{l+1} ( \eufrak{s}_{l+1} ) 
 \right] 
\\ 
\nonumber 
 & & 
 \! \! \! \! \! \! \! 
 + 
 \sum_{d=1}^l 
 \sum_{0 = k_d < \cdots < k_0 = n} 
 c_{ \eufrak{K}_d } 
 E_\sigma \left[  \int_{X^d} 
 \varepsilon^+_{\eufrak{s}_d} 
 G 
 \prod_{p=1}^d 
 \varepsilon^+_{\eufrak{s}_d \setminus s_p} 
 u_{s_p}^{k_{p-1} - k_p } 
 \hskip0.05cm 
 d\sigma^d ( \eufrak{s}_d ) 
 \right] 
\\ 
\nonumber 
 & & 
 \! \! \! \! \! \! \! 
 - 
 \sum_{d=1}^{l+1} 
 \sum_{0 = k_d < \cdots < k_0 = n} 
 \! \! \! \! \! \! \! 
 c_{ \eufrak{K}_d } 
 E_\sigma \left[  \int_{X^d} 
 \delta_\sigma ( 
 \varepsilon^+_{ \eufrak{s}_{d-1}} 
 u 
 )^{k_{d-1}-1}
 \varepsilon^+_{ \eufrak{s}_{d-1}} 
 ( 
 u_{s_d} 
 G 
 ) 
 \prod_{p=1}^{d-1} 
 \varepsilon^+_{ \eufrak{s}_{d-1} \setminus s_p} 
 u_{s_p}^{k_{p-1} - k_p} 
 \hskip0.05cm 
 d\sigma^d ( \eufrak{s}_d ) 
 \right] 
\\ 
 & = & 
 {\cal A}_l + \sum_{d=1}^l {\cal B}_d - \sum_{d=1}^{l+1} {\cal C}_d
, 
\end{eqnarray} 
 which will be proved by induction on $l = 0 , \ldots , n$. 
 First, note that \eqref{1n2} holds for $l=0$ 
 as by \eqref{prr} and 
 \eqref{dr1} we have 
\begin{eqnarray*} 
\lefteqn{ 
 E_\sigma \left[ 
 \delta_\sigma ( u )^n 
 G 
 \right] 
 = 
 E_\sigma \left[ 
 \int_X 
 u_{s_1} 
 D_{s_1} 
 ( 
 \delta_\sigma ( u)^{n-1} 
 G 
 ) 
 \sigma ( ds_1 ) 
 \right] 
} 
\\ 
 & = & 
 E_\sigma \left[ 
 \int_X 
 u_{s_1} 
 \varepsilon^+_{s_1} 
 \delta_\sigma ( u)^{n-1} 
 \varepsilon^+_{s_1} 
 G 
 \sigma ( ds_1 ) 
 \right] 
 - 
 E_\sigma \left[ 
 G 
 \int_X 
 u_{s_1} 
 \delta_\sigma ( u )^{n-1} 
 \sigma ( ds_1 ) 
 \right] 
\\ 
 & = & 
 E_\sigma \left[ 
 \int_X 
 u_{s_1} 
 ( 
 u_{s_1} 
 + 
 \delta_\sigma 
 ( 
 \varepsilon^+_{s_1} u 
 ) 
 )^{n-1} 
 \varepsilon^+_{s_1} 
 G 
 \sigma ( ds_1 ) 
 \right]  
 - 
 E_\sigma \left[ 
 G 
 \int_X 
 u_{s_1} 
 \delta_\sigma ( u)^{n-1} 
 \sigma ( ds_1 ) 
 \right] 
\\ 
 & = & 
 {\cal A}_0 - {\cal C}_1
, 
\end{eqnarray*} 
 which also proves the lemma in case $n=1$. 
 Next, when $n\geq 2$, for $l = 0 , \ldots , n-1$, 
 using the duality formula \eqref{dr1} and the relations 
$ 
 \varepsilon^+_{s_{l+1}} 
 \delta_\sigma ( 
 \varepsilon^+_{ \eufrak{s}_l} 
 u 
 ) 
 = 
 \varepsilon^+_{ \eufrak{s}_l} 
 u_{s_{l+1}} 
 + 
 \delta_\sigma ( 
 \varepsilon^+_{ \eufrak{s}_{l+1}} 
 u 
 ) 
$, 
 cf. \eqref{commrel}, 
 and $D_{s_{l+2}} = \varepsilon^+_{s_{l+2}} - I$, we rewrite the 
 first term in \eqref{1n2} as 
\begin{eqnarray*} 
\lefteqn{ 
 {\cal A}_l = \sum_{0 = k_{l+1} < \cdots < k_0 = n} 
 c_{ \eufrak{K}_{l+1} } 
} 
\\ 
 & & 
 \! \! \! \! \! \! \! \! \! 
 E_\sigma \left[ 
 \int_{X^{l+1}} 
 \varepsilon^+_{ \eufrak{s}_{l+1} } 
 G 
 \varepsilon^+_{ \eufrak{s}_l } 
 u_{s_{l+1}} 
 \left( 
 \varepsilon^+_{ \eufrak{s}_l} 
 u_{s_{l+1}}  
 + 
 \delta_\sigma ( 
 \varepsilon^+_{ \eufrak{s}_{l+1}} 
 u 
 ) 
 \right)^{k_l-1}
 \prod_{p=1}^l 
 \varepsilon^+_{ \eufrak{s}_{l+1} \setminus s_p} 
 u_{s_p}^{k_{p-1} - k_p} 
 \hskip0.05cm 
 d\sigma^{l+1} ( \eufrak{s}_{l+1} ) 
 \right] 
\\ 
 & = & 
 \! \! \! \! \! \! \! \! \! 
 \sum_{0 \leq k_{l+1} < \cdots < k_0 = n} 
 \! \! \! \! \! \! \! 
 c_{ \eufrak{K}_{l+1} } 
 E_\sigma \left[  \int_{X^{l+1}} 
 \delta_\sigma ( 
 \varepsilon^+_{ \eufrak{s}_{l+1}} 
 u 
 )^{k_{l+1}}
 \varepsilon^+_{ \eufrak{s}_{l+1}} 
 G 
 \varepsilon^+_{ \eufrak{s}_l} 
 u_{s_{l+1}}^{k_l-k_{l+1}} 
 \prod_{p=1}^l 
 \varepsilon^+_{ \eufrak{s}_{l+1} \setminus s_p} 
 u_{s_p}^{k_{p-1} - k_p} 
 \hskip0.05cm 
 d\sigma^{l+1} ( \eufrak{s}_{l+1} ) 
 \right] 
\\ 
 & = & 
 \! \! \! \! \! \! \! \! \! 
 \sum_{1 \leq k_{l+1} < \cdots < k_0 = n} 
 \! \! \! \! \! \! \! \! \! 
 c_{ \eufrak{K}_{l+1} } 
 E_\sigma \left[  \int_{X^{l+1}} 
 \! \! \! \! 
 \delta_\sigma ( 
 \varepsilon^+_{ \eufrak{s}_{l+1}} 
 u 
 )^{k_{l+1}}
 \varepsilon^+_{ \eufrak{s}_{l+1}} 
 G 
 \prod_{p=1}^{l+1} 
 \varepsilon^+_{ \eufrak{s}_{l+1} \setminus s_p} 
 u_{s_p}^{k_{p-1} - k_p} 
 \hskip0.05cm 
 d\sigma^{l+1} ( \eufrak{s}_{l+1} ) 
 \right] 
\\ 
 & & 
 + 
 \sum_{0=k_{l+1} < \cdots < k_0 = n} 
 c_{ \eufrak{K}_{l+1} } 
 E_\sigma \left[  \int_{X^{l+1}} 
 \varepsilon^+_{ \eufrak{s}_{l+1} } 
 G 
 \prod_{p=1}^{l+1} 
 \varepsilon^+_{ \eufrak{s}_{l+1} \setminus s_p } 
 u_{s_p}^{k_{p-1} - k_p} 
 \hskip0.05cm 
 d\sigma^{l+1} ( \eufrak{s}_{l+1} ) 
 \right] 
\\ 
 & = & 
 \sum_{1 \leq k_{l+1} < \cdots < k_0 = n} 
 c_{ \eufrak{K}_{l+1} } 
\\ 
 & & 
 \! \! \! \! \! \! \! \! \! \! 
  E_\sigma \left[  \int_{X^{l+2}} 
 \varepsilon^+_{ \eufrak{s}_{l+1}} 
 u_{s_{l+2}} 
 D_{s_{l+2}} 
 \left( 
 \delta_\sigma ( 
 \varepsilon^+_{ \eufrak{s}_{l+1}} 
 u 
 )^{k_{l+1}-1}
 \varepsilon^+_{ \eufrak{s}_{l+1}} 
 G 
 \prod_{p=1}^{l+1} 
 \varepsilon^+_{ \eufrak{s}_{l+1} \setminus s_p} 
 u_{s_p}^{k_{p-1} - k_p} 
 \right) 
 \hskip0.05cm 
 d\sigma^{l+2} ( \eufrak{s}_{l+2} ) 
 \right] 
\\ 
 & & 
 + 
 \sum_{0=k_{l+1} < \cdots < k_0 = n} 
 c_{ \eufrak{K}_{l+1} } 
 E_\sigma \left[  \int_{X^{l+1}} 
 \varepsilon^+_{ \eufrak{s}_{l+1}} 
 G 
 \prod_{p=1}^{l+1} 
 \varepsilon^+_{ \eufrak{s}_{l+1} \setminus s_p} 
 u_{s_p}^{k_{p-1} - k_p} 
 d\sigma^{l+1} ( \eufrak{s}_{l+1} ) 
 \right] 
\\ 
 & = & 
 \sum_{1 \leq k_{l+1} < \cdots < k_0 = n} 
 c_{ \eufrak{K}_{l+1} } 
\\ 
 & & 
 \! \! \! \! \! \! \! \! \! \! \! \! 
 E_\sigma \left[  \int_{X^{l+2}} 
 \varepsilon^+_{ \eufrak{s}_{l+1}} 
 u_{s_{l+2}} 
 \right. 
 \left. 
 \varepsilon^+_{s_{l+2} } 
 \left( 
 \delta_\sigma ( 
 \varepsilon^+_{ \eufrak{s}_{l+1}} 
 u 
 )^{k_{l+1}-1}
 \varepsilon^+_{ \eufrak{s}_{l+1}} 
 G 
 \prod_{p=1}^{l+1} 
 \varepsilon^+_{ \eufrak{s}_{l+1} \setminus s_p} 
 u_{s_p}^{k_{p-1} - k_p} 
 \right) 
 d\sigma^{l+2} ( \eufrak{s}_{l+2} ) 
 \right] 
\\ 
 & & 
 \! \! \! \! \! \! \! \! \! 
 - 
 \sum_{1 \leq k_{l+1} < \cdots < k_0 = n} 
 c_{ \eufrak{K}_{l+1} } 
 E_\sigma \left[  \int_{X^{l+2}} 
 \delta_\sigma ( 
 \varepsilon^+_{ \eufrak{s}_{l+1}} 
 u 
 )^{k_{l+1}-1}
 \varepsilon^+_{ \eufrak{s}_{l+1}} 
 ( 
 u_{s_{l+2}} 
 G 
 ) 
 \prod_{p=1}^{l+1} 
 \varepsilon^+_{ \eufrak{s}_{l+1} \setminus s_p} 
 u_{s_p}^{k_{p-1} - k_p} 
 \hskip0.05cm 
 d\sigma^{l+2} ( \eufrak{s}_{l+2} ) 
 \right] 
\\ 
 & & 
 \! \! \! \! \! \! \! \! \! 
 + 
 \sum_{0 = k_{l+1} < \cdots < k_0 = n} 
 c_{ \eufrak{K}_{l+1} } 
 E_\sigma \left[  \int_{X^{l+1}} 
 \varepsilon^+_{ \eufrak{s}_{l+1}} 
 G 
\right. 
\left. 
 \prod_{p=1}^{l+1} 
 \varepsilon^+_{ \eufrak{s}_{l+1} \setminus s_p} 
 u_{s_p}^{k_{p-1} - k_p} 
 \hskip0.05cm 
 d\sigma^{l+1} ( \eufrak{s}_{l+1} ) 
 \right] 
\\ 
 & = & 
 {\cal A}_{l+1} 
 + 
 {\cal B}_{l+1} 
 - 
 {\cal C}_{l+2} 
, 
\end{eqnarray*} 
 which proves \eqref{1n2} by induction on $l = 1,\ldots , n-1$, as 
$$ 
 E_\sigma \left[ 
 \delta_\sigma ( u )^n 
 G 
 \right] 
 = 
 {\cal A}_0 - {\cal C}_1 
 = 
 - 
 {\cal C}_1 
 + 
 \sum_{d=0}^{n-1} {\cal A}_d - {\cal A}_{d+1} 
 = 
 - 
 {\cal C}_1 
 + 
 \sum_{d=0}^{n-1} {\cal B}_{d+1} - {\cal C}_{d+2} 
 = 
 \sum_{d=1}^n {\cal B}_d 
 - 
 \sum_{d=1}^{n+1} {\cal C}_d. 
$$ 
\end{Proof} 
\begin{Proofy} {\em \hskip-0.12cm of Theorem~\ref{l22}}. 
 We check that in \eqref{plkmn}, all terms with $a=0$ and $0 \leq b \leq n-1$ 
 vanish, hence in particular the formula also holds when $n=0$. 
 When $n\geq 1$ the proof of \eqref{plkmn} is obtained by application to $c=n$ 
 or $c=n+1$ of the following identity: 
\begin{eqnarray} 
\label{den} 
\lefteqn{ 
 E_\sigma \left[ 
 \delta_\sigma ( u )^n 
 F 
 \right] 
 = 
 (-1)^c 
 \sum_{a=0}^{n-c} 
 \sum_{l_1 + \cdots + l_{a+1} = n -c - a 
 \atop l_1 ,\ldots ,l_{a+1} \geq 0} 
 C_{\eufrak{L}_{a+1},a+c} 
} 
\\ 
\nonumber 
 & & 
 \ \ \ \ \ \ \ \ \ \ \ \ \ \ \ \ \ 
 E_\sigma \left[ 
 \int_{X^{a+c}} 
 \delta_\sigma ( 
 \varepsilon^+_{ \eufrak{s}_a} 
 u 
 )^{l_{a+1}}
 \varepsilon^+_{ \eufrak{s}_a} 
 F 
 \prod_{q=a+1}^{a+c} 
 \varepsilon^+_{ \eufrak{s}_a} 
 u_{s_q} 
 \prod_{p=1}^a 
 \varepsilon^+_{ \eufrak{s}_a \setminus s_p} 
 u_{s_p}^{1+l_p} 
 \hskip0.05cm 
 d\sigma^{a+c} ( \eufrak{s}_{a+c} ) 
 \right] 
\\ 
\nonumber 
 & & 
 \! \! \! \! \! \! \! \! \! \! \! \! 
 + 
 \sum_{b=0}^{c-1} 
 (-1)^b 
 \sum_{a=0}^{n-b} 
 \sum_{l_1 + \cdots + l_a = n - b - a 
 \atop l_1,\ldots ,l_a \geq 0} 
 C_{\eufrak{L}_a,a+b} 
 E_\sigma \left[  \int_{X^{a+b}} 
 \varepsilon^+_{ \eufrak{s}_a} 
 F 
 \prod_{q=a+1}^{a+b} 
 \varepsilon^+_{ \eufrak{s}_a} 
 u_{s_q} 
 \prod_{p=1}^a 
 \varepsilon^+_{ \eufrak{s}_a \setminus s_p} 
 u_{s_p}^{1+l_p} 
 \hskip0.05cm 
 d\sigma^{a+b} ( \eufrak{s}_{a+b} ) 
 \right] 
\\ 
\nonumber 
 & = & 
 {\cal D}_c + \sum_{b=0}^{c-1} {\cal E}_b
, 
\end{eqnarray} 
 which will be proved by induction on $c=1, \ldots ,n+1$. 
 First, we note that since 
$$ 
 C_{\eufrak{L}_a,a} 
 = 
 \prod_{p=1}^a 
 {l_1 + \cdots + l_p + p - 1  
 \choose l_1 + \cdots + l_{p-1} + p - 1 } 
, 
$$ 
 the identity \eqref{den} holds for $c=1$ from Lemma~\ref{l11}. 
 Next, for all $c=1,\ldots ,n-1$, 
 applying Lemma~\ref{l11} with $n=l_{a+1}$ and 
$$ 
 G = 
 \varepsilon^+_{ \eufrak{s}_a} 
 F 
 \prod_{q=a+1}^{a+c} 
 \varepsilon^+_{ \eufrak{s}_a} 
 u_{s_q} 
 \prod_{p=1}^a 
 \varepsilon^+_{ \eufrak{s}_a \setminus s_p} 
 u_{s_p}^{1+l_p} 
$$ 
 and fixing $s_1 ,\ldots,s_{a+c}$, 
 we rewrite the first term in \eqref{den} using \eqref{defin} 
 and the change of index 
$$ 
 k_{d-p} = p + m_1 + \cdots + m_p, \qquad 0 \leq p \leq d
, 
$$ 
 as 
\begin{eqnarray} 
\nonumber 
\lefteqn{ 
 {\cal D}_c = (-1)^c 
 \sum_{a=0}^{n-c} 
 \sum_{l_1 + \cdots + l_{a+1} = n - c - a 
 \atop l_1 ,\ldots ,l_{a+1} \geq 0} 
 C_{\eufrak{L}_{a+1},a+c} 
} 
\\ 
\nonumber 
 & & 
 \times 
 E_\sigma \left[  \int_{X^{a+c}} 
 \delta_\sigma ( 
 \varepsilon^+_{ \eufrak{s}_a} 
 u 
 )^{l_{a+1}}
 \varepsilon^+_{ \eufrak{s}_a} 
 F 
 \prod_{q=a+1}^{a+c} 
 \varepsilon^+_{ \eufrak{s}_a} 
 u_{s_q} 
 \prod_{p=1}^a 
 \varepsilon^+_{ \eufrak{s}_a \setminus s_p} 
 u_{s_p}^{1+l_p} 
 \hskip0.05cm 
 d\sigma^{a+c} ( \eufrak{s}_{a+c} ) 
 \right] 
\\ 
\nonumber 
 & = & 
 (-1)^c 
 \sum_{a=0}^{n-c} 
 \sum_{l_1 + \cdots + l_{a+1} = n - c - a  
 \atop l_1 ,\ldots , l_{a+1} \geq 0} 
 C_{\eufrak{L}_{a+1},a+c} 
 \sum_{d=0}^{l_{a+1}} 
 \sum_{m_1 + \cdots + m_d = l_{a+1} - d 
 \atop m_1,\ldots ,m_d \geq 0} 
 \prod_{p=1}^d 
 {m_1 + \cdots + m_p + p - 1  
 \choose m_1 + \cdots + m_{p-1} + p - 1 } 
\\ 
\nonumber 
 & & 
 \times E_\sigma \left[  \int_{X^{a+c+d}} 
 \varepsilon^+_{ \eufrak{s}_{a+d}} 
 F 
 \prod_{q=a+d+1}^{a+c+d} 
 \varepsilon^+_{ \eufrak{s}_{a+d}} 
 u_{s_q} 
\right. 
 \left. 
 \prod_{p=1}^a 
 \varepsilon^+_{ \eufrak{s}_{a+d} \setminus s_p } 
 u_{s_p}^{1+l_p} 
 \prod_{k=a+1}^{a+d} 
 \varepsilon^+_{ \eufrak{s}_{a+d} \setminus s_p } 
 u_{s_k}^{1+m_{k-a}} 
 \hskip0.05cm 
 d\sigma^{a+c+d} ( \eufrak{s}_{a+c+d} ) 
 \right] 
\\ 
\nonumber 
 & & 
 - 
 (-1)^c 
 \sum_{a=0}^{n-c} 
 \sum_{l_1 + \cdots + l_{a+1} = n - c - a  
 \atop l_1 ,\ldots ,l_{a+1} \geq 0} 
 C_{\eufrak{L}_{a+1},a+c} 
 \sum_{d=1}^{l_{a+1}} 
 \sum_{m_1 + \cdots + m_d = l_{a+1} - d 
 \atop m_1,\ldots ,m_d \geq 0} 
 \prod_{p=1}^d 
 {m_1 + \cdots + m_p + p - 1  
 \choose m_1 + \cdots + m_{p-1} + p - 1} 
\\ 
\nonumber 
 & & 
 \times E_\sigma \left[  \int_{X^{a+c+d}} 
 \delta_\sigma ( 
 \varepsilon^+_{ \eufrak{s}_{a+d-1} } 
 u 
 )^{m_d}
 \varepsilon^+_{ \eufrak{s}_{a+d-1} } 
 F 
 \prod_{q=a+d}^{a+c+d} 
 \varepsilon^+_{ \eufrak{s}_{a+d-1} } 
 u_{s_q} 
\right. 
\\ 
\nonumber 
 & & 
\left. 
 \ \ \ \ \  \ \ \ \ \  \ \ \ \ \  \ \ \ \ \  \ \ \ \ \  \ \ \ \ \ 
 \times 
 \prod_{p=1}^a 
 \varepsilon^+_{ \eufrak{s}_{a+d-1} \setminus s_p } 
 u_{s_p}^{1+l_p} 
 \prod_{k=a+c+1}^{a+c+d-1} 
 \varepsilon^+_{ \eufrak{s}_{a+d-1} \setminus s_p} 
 u_{s_k}^{1+m_{k-a-c}} 
 \hskip0.05cm 
 d\sigma^{a+c+d} ( \eufrak{s}_{a+c+d} ) 
 \right] 
\\ 
\nonumber 
 & = & 
 (-1)^c 
 \sum_{a'=0}^{n-c} 
 \sum_{l'_1 + \cdots + l'_{a'} = n - c - a' 
 \atop l'_1 ,\ldots ,l'_{a'} \geq 0} 
 \! \! \! \! \! \! \! 
 C_{\eufrak{L}_{a'},a'+c} 
 E_\sigma \left[ 
 \int_{X^{a'+c}} 
 \! \! 
 \varepsilon^+_{ \eufrak{s}_{a'}} 
 F 
 \prod_{q=a'+1}^{a'+c} 
 \varepsilon^+_{ \eufrak{s}_{a'+c}} 
 u_{s_q} 
 \prod_{p=1}^{a'} 
 \varepsilon^+_{ \eufrak{s}_{a'} \setminus s_p} 
 u_{s_p}^{1+l'_p} 
 \hskip0.05cm 
 d\sigma^{a'+c} ( \eufrak{s}_{a'+c} ) 
 \right] 
\\ 
\label{113} 
\\ 
\label{223} 
 & & 
 + 
 (-1)^{c+1} 
 \sum_{a'=0}^{n-c} 
 \sum_{l'_1 + \cdots + l'_{a'+1} = n - c - a' - 1 
 \atop l'_1 , \ldots ,l'_{a'+1} \geq 0} 
 C_{\eufrak{L}_{a'+1},a'+c+1} 
\\ 
\nonumber 
 & & 
 \ \ \ \ \  \ \ \ 
 \times 
 E_\sigma \left[ 
 \int_{X^{a'+c+1}} 
 \delta_\sigma ( 
 \varepsilon^+_{ \eufrak{s}_{a'}} 
 u 
 )^{l'_{a'+1}}
 \varepsilon^+_{ \eufrak{s}_{a'}} 
 F 
 \prod_{q=a'+1}^{a'+c+1} 
 \varepsilon^+_{ \eufrak{s}_{a'}} 
 u_{s_q} 
 \prod_{p=1}^{a'} 
 \varepsilon^+_{ \eufrak{s}_{a'} \setminus s_p} 
 u_{s_p}^{1+l'_p} 
 \hskip0.05cm 
 d\sigma^{a'+c+1} ( \eufrak{s}_{a'+c+1} ) 
 \right] 
\\ 
\nonumber 
 & = & 
 {\cal E}_c + {\cal D}_{c+1}, 
\end{eqnarray} 
 under the changes of indices 
$$ 
 l'_1+\cdots + l'_{a'} = l_1+\cdots + l_a 
 + m_1+\cdots + m_d 
, 
 \qquad 
 a'=a+d, 
$$ 
 in \eqref{113} when $d=0,\ldots,l_{a+1}$, and 
$$ 
 l'_1+\cdots + l'_{a'+1} = l_1+\cdots + l_a + m_1+\cdots + m_d, 
 \qquad 
 a'+1=a+d, 
$$ 
 in \eqref{223} when $d = 1 , \ldots , l_{a+1}$. 
 Noting that in \eqref{223}, 
 the summation on $a'$ actually ends at $a' = n-c-1$ when $c<n$. 
 We conclude the proof by induction, as 
$$ 
 E_\sigma \left[ 
 \delta_\sigma ( u )^n 
 F 
 \right] 
 = 
 {\cal D}_1 - {\cal D}_{n+1} + {\cal E}_0 
 = 
 {\cal E}_0 + 
 \sum_{b=1}^n {\cal D}_b - {\cal D}_{b+1} 
 = 
 \sum_{b=0}^n {\cal E}_b 
, 
$$ 
 and by the change of indices $(a,b)\to (a,b-a)$ in \eqref{plkmn}. 
\end{Proofy} 
\section{Recursive moment identities} 
\label{rmi} 
 The main results of this section are Propositions~\ref{p01} 
 and \ref{c1}. 
 Their proofs are stated using Lemma~\ref{l5} above 
 and Proposition~\ref{pr1} below, 
 and they are used to prove the main results of 
 Section~\ref{mr}. 
 In the next theorem we use the notation $\Delta_s$ of Definition~\ref{d1} 
 and let 
$$ 
 \Delta_{\eufrak{s}_{j} } 
 = 
 \Delta_{\eufrak{s}_0 } 
 \cdots 
 \Delta_{\eufrak{s}_{j} } 
, 
 \qquad 
 \eufrak{s}_{j} = ( {s}_0 , \ldots , {s}_{j} ), 
$$ 
 and 
$$ 
 d\sigma^{b+1} ( \eufrak{s}_b ) 
 = 
 \sigma ( d {s}_0 ) 
 \cdots 
 \sigma ( d {s}_b ) 
, 
 \qquad 
 \eufrak{s}_{b} = ( {s}_0 , \ldots , {s}_{b} ), 
$$ 
 $0 \leq j \leq b$. 
\begin{prop} 
\label{p01} 
 Let $N\geq 0$ and let $u\in \lee_{2,1}$ 
 be bounded with 
 $ 
 \displaystyle 
 u \in \bigcap_{p=1}^{N+1} 
 L^\infty ( \Omega^X , L^p_\sigma ( X ) )$ 
 and 
$$ 
 E_\sigma 
 \left[ 
 \int_{X^{b+1}} 
 \Big| 
 \Delta_{s_0} 
 \cdots  
 \Delta_{s_j} 
 \left( 
 \prod_{q=a+1}^b 
 u_{s_q} 
 \prod_{p=0}^a 
 u^{l_p}_{s_p} 
 \right) 
 \Big| 
 d\sigma^{b+1} ( \eufrak{s}_b ) 
 \right] 
 < \infty 
, 
$$ 
 $l_0 + \cdots + l_a \leq N+1$, 
 $l_0 ,\ldots ,l_a \geq 1$, 
 $0 \leq j \leq a \leq b \leq N$. 
 Then for all $n=0, \ldots ,N$ we have 
 $\delta_\sigma ( u ) \in L^{n+1} ( \Omega^X , \pi_\sigma )$ 
 and 
\begin{eqnarray} 
\label{xt} 
\lefteqn{ 
 E_\sigma [ \delta_\sigma ( u )^{n+1} ] 
 = 
 \sum_{k=0}^{n-1} 
 {n \choose k} 
 E_\sigma \left[ 
 \delta_\sigma ( u )^k 
 \int_X 
 u^{n-k+1}_t
 \sigma ( dt ) 
 \right] 
} 
\\ 
\nonumber 
 & & 
\! \! \! \! \! \! \! \! 
 + 
 \sum_{a=0}^n 
 \sum_{j=0}^a 
 \sum_{b = a}^n 
 \sum_{{l_0 + \cdots + l_a = n - b 
 \atop l_0 ,\ldots ,l_a \geq 0 
 } 
 } 
 {a \choose j} 
 C_{\eufrak{L}_a,b}^{l_0,n} 
 E_\sigma 
 \left[ 
 \int_{X^{b+1}} 
 \Delta_{\eufrak{s}_{j} } 
 \left( 
 \prod_{q=a+1}^b 
 u_{s_q} 
 \prod_{p=0}^a 
 u^{1+l_p}_{s_p} 
 \right) 
 d\sigma^{b+1} ( \eufrak{s}_b ) 
 \right] 
, 
\end{eqnarray} 
 where 
$$ 
 C_{\eufrak{L}_a,b}^{l_0,n} 
 = 
 (-1)^{b-a} 
 {n \choose l_0}  
 C_{\eufrak{L}_a , b } 
, 
$$ 
 and $C_{\eufrak{L}_a , b }$ is defined in \eqref{defin}. 
\end{prop} 
\begin{Proof} 
 When $u:\Omega^X \times X \to \real$ is a 
 bounded process with compact support in $X$ 
 this result is a direct consequence 
 of Lemma~\ref{l5} and Proposition~\ref{pr1} below applied with 
 $n=k$ and $l_0 + k - b = n - b$. 
 We conclude the proof by induction and a limiting argument, 
 as follows. 
 Let $(K_r)_{r\geq 1}$ denote an increasing family of compact subsets of $X$ 
 such that $X = \displaystyle \bigcup_{r\geq 1} K_r$. 
 The family of processes 
 $u^{(r)}_x (\omega ) : = 
 u_x 
 {\bf 1}_{K_r} (x) 
 $, 
 $r \geq 1$, converges in $\lee_{2,1}$ to $u$ as $r$ goes 
 to infinity, hence 
 $\delta_\sigma (u^{(r)})$ converges to $\delta (u)$ 
 in $L^2 ( \Omega^X , \pi_\sigma )$ as $r$ goes to infinity. 
 Clearly the result holds for $N=0$ 
 by applying the formula to the process 
 $u^{(r)}$ which is bounded with compact support 
 by letting $r$ go to infinity. 
 Next, letting $N\geq 0$ and assuming 
 that 
 $\delta_\sigma ( u ) \in L^{n+1} ( \Omega^X , \pi_\sigma )$ 
 and that \eqref{xt} holds for all $n=0, \ldots ,N$, 
 we note that for all even integer $m \in \{ 
 2, \ldots , N+1\}$ we have the bound 
\begin{eqnarray*} 
\lefteqn{ 
 E_\sigma [ \delta_\sigma ( u )^m ] 
 \leq 
 \sum_{k=0}^{m-2} 
 {m-1 \choose k} 
 E_\sigma \left[ 
 \delta_\sigma ( u )^{m-2} 
 \right]^{k/(m-2)} 
 \left\| 
 \int_X 
 | u_t |^{m-k} 
 \sigma ( dt ) 
 \right\|_\infty 
} 
\\ 
 & & 
 \!  \!  \!  \!  \!  \!  \!  \!  \!  \!  \! 
 + 
 \sum_{a=0}^{m-1} 
 \sum_{j=0}^a 
 \sum_{b = a}^{m-1} 
 \sum_{{l_0 + \cdots + l_a = m - b - 1 
 \atop l_0 ,\ldots ,l_a \geq 0 
 } 
 } 
 \!  \!  \! 
 {a \choose j} 
 C_{\eufrak{L}_a,b}^{l_0,m-1} 
 E_\sigma 
 \left[ 
 \int_{X^{b+1}} 
 \left| 
 \Delta_{\eufrak{s}_{j} } 
 \left( 
 \prod_{q=a+1}^b 
 u_{s_q} 
 \prod_{p=0}^a 
 u^{1+l_p}_{s_p} 
 \right) 
 \right| 
 d\sigma^{b+1} ( \eufrak{s}_b ) 
 \right] 
, 
\end{eqnarray*} 
 which, applied to $u^{(r)} (\omega )$, allows us to extend 
 \eqref{xt} to the order $N+1$ by uniform integrability 
 after taking the limit as $r$ goes to infinity. 
\end{Proof} 
 Let us consider some particular cases of Proposition~\ref{p01}. 
 For $n=1$, Relation~\eqref{xt} reads 
\small 
\begin{align*} 
 E_\sigma \left[ 
 \delta_\sigma ( u )^2 
 \right] 
 = \ & 
 E_\sigma \left[  \int_X 
 | u_s |^2 
 \sigma ( ds ) 
 \right] 
 & 
\\ 
 & 
 - 
 E_\sigma \left[  \int_{X^2} 
 \Delta_{s_1} 
 ( 
 u_{s_1}  
 u_{s_2} 
 ) 
 \sigma ( ds_1 ) \sigma ( ds_2 ) 
 \right] 
 & [a=0,b=1,j=0] 
\\ 
 & 
 + 
 E_\sigma \left[  \int_{X^2} 
 \Delta_{s_1} 
 ( 
 u_{s_1}  
 u_{s_2} 
 ) 
 \sigma ( ds_1 ) \sigma ( ds_2 ) 
 \right] 
 & [a=1,b=1,j=0] 
\\ 
 & 
 + 
 E_\sigma \left[  \int_{X^2} 
 \Delta_{s_1} 
 \Delta_{s_2} 
 ( 
 u_{s_1}  
 u_{s_2} 
 ) 
 \sigma ( ds_1 ) \sigma ( ds_2 ) 
 \right] 
 & [a=1,b=1,j=1] 
\\ 
 = \ & 
 E_\sigma \left[  \int_X 
 | u_s |^2 
 \sigma ( ds ) 
 \right] 
 + 
 E_\sigma \left[  \int_{X^2} 
 \Delta_{s_1} 
 \Delta_{s_2} 
 ( 
 u_{s_1}  
 u_{s_2} 
 ) 
 \sigma ( ds_1 ) \sigma ( ds_2 ) 
 \right] 
, 
 & 
\end{align*} 
\normalsize 
\baselineskip0.7cm 
 which coincides with \eqref{qas1}. 
 On the other hand for $n=2$ Relation~\eqref{xt} yields 
 the third moment 
\small 
\begin{align} 
\nonumber 
 & 
 E_\sigma \left[ 
 \delta_\sigma ( u )^3 
 \right] 
 = 
 E_\sigma \left[ 
 \int_X 
 u^3_s 
 \sigma ( ds ) 
 \right] 
 + 
 2 
 E_\sigma \left[ 
 \delta ( u ) 
 \int_{X} 
 u_s^2 
 \sigma ( ds ) 
 \right] 
\\ 
\nonumber 
 & 
 - 
 2 
 E_\sigma \left[  \int_{X^2} 
 \Delta_{s_0} 
 ( 
 u^2_{s_0} 
 u_{s_1} 
 ) 
 \sigma ( ds_0 ) \sigma ( ds_1 ) 
 \right] 
 & [a=0,b=1,j=0] 
\\ 
\nonumber 
 & 
 + 
 E_\sigma \left[  \int_{X^2} 
 \Delta_{s_0} 
 ( 
 u_{s_0} 
 u_{s_1} 
 u_{s_2} 
 ) 
 \sigma ( ds_0 ) \sigma ( ds_1 ) 
 \right] 
 & [a=0,b=2,j=0] 
\\ 
\nonumber 
 & 
 + 
 E_\sigma \left[  \int_{X^2} 
 \Delta_{s_0} 
 ( 
 u_{s_0} 
 u_{s_1}^2 
 ) 
 \sigma ( ds_0 ) \sigma ( ds_1 ) 
 \right] 
 + 
 2 
 E_\sigma \left[  \int_{X^2} 
 \Delta_{s_0} 
 ( 
 u_{s_0}^2 
 u_{s_1} 
 ) 
 \sigma ( ds_0 ) \sigma ( ds_1 ) 
 \right] 
 & [a=1,b=1,j=0] 
\\ 
\nonumber 
 & 
 + 
 3 
 E_\sigma \left[  \int_{X^2} 
 \Delta_{s_0} \Delta_{s_1} 
 ( 
 u_{s_0} 
 u_{s_0}^2 
 ) 
 \sigma ( ds_0 ) \sigma ( ds_1 ) 
 \right] 
 & [a=1,b=1,j=1] 
\\ 
\nonumber 
 & 
 - 
 E_\sigma \left[  \int_{X^2} 
 \Delta_{s_0} 
 ( 
 u_{s_0} 
 u_{s_1} 
 u_{s_2} 
 ) 
 \sigma ( ds_0 ) \sigma ( ds_1 ) 
 \right] 
 & [a=1,b=2,j=0] 
\\ 
\nonumber 
 & 
 - 
 E_\sigma \left[  \int_{X^3} 
 \Delta_{s_0} \Delta_{s_1} 
 ( 
 u_{s_0} 
 u_{s_1} 
 u_{s_2} 
 ) 
 \sigma ( ds_0 ) \sigma ( ds_1 ) \sigma ( ds_2 ) 
 \right] 
 & [a=1,b=2,j=1] 
\\ 
\nonumber 
 & 
 + 
 E_\sigma \left[  \int_{X^3} 
 \Delta_{s_0} 
 ( 
 u_{s_0} 
 u_{s_1} 
 u_{s_2} 
 ) 
 \sigma ( ds_0 ) \sigma ( ds_1 ) \sigma ( ds_2 ) 
 \right] 
 & [a=2,b=2,j=0] 
\\ 
\nonumber 
 & 
 + 
 E_\sigma \left[  \int_{X^3} 
 \Delta_{s_0} \Delta_{s_1} 
 ( 
 u_{s_0} 
 u_{s_1} 
 u_{s_2} 
 ) 
 \sigma ( ds_0 ) \sigma ( ds_1 ) \sigma ( ds_2 ) 
 \right] 
 & [a=2,b=2,j=1] 
\\ 
\nonumber 
 & 
 + 
 E_\sigma \left[  \int_{X^3} 
 \Delta_{s_0} \Delta_{s_1} \Delta_{s_2} 
 ( 
 u_{s_0} 
 u_{s_1} 
 u_{s_2} 
 ) 
 \sigma ( ds_0 ) \sigma ( ds_1 ) \sigma ( ds_2 ) 
 \right] 
 & [a=2,b=2,j=2] 
\end{align} 
\vskip-0.6cm 
\begin{align} 
\nonumber 
 = & \ 
 E_\sigma \left[ 
 \int_X 
 u^3_s 
 \sigma ( ds ) 
 \right] 
 + 
 2 
 E_\sigma \left[ 
 \delta ( u ) 
 \int_{X} 
 u_s^2 
 \sigma ( ds ) 
 \right] 
 + 
 E_\sigma \left[  \int_{X^2} 
 u_{s_0} 
 D_{s_0} 
 u_{s_1}^2 
 \sigma ( ds_0 ) \sigma ( ds_1 ) 
 \right] 
 & 
\\ 
\label{qas3} 
 & 
 + 
 3 
 E_\sigma \left[  \int_{X^2} 
 \Delta_{s_0} \Delta_{s_1} 
 ( 
 u_{s_0} 
 u_{s_1}^2 
 ) 
 \sigma ( ds_0 ) \sigma ( ds_1 ) 
 \right] 
 + 
 E_\sigma \left[  \int_{X^3} 
 \Delta_{s_0} \Delta_{s_1} \Delta_{s_2} 
 ( 
 u_{s_0} 
 u_{s_1} 
 u_{s_2} 
 ) 
 \sigma ( ds_0 ) \sigma ( ds_1 ) \sigma ( ds_2 ) 
 \right] 
, 
 & 
\end{align} 
\normalsize 
\baselineskip0.7cm 
 which recovers \eqref{qas2} by the duality relation \eqref{dr1}. 
 As a consequence of Proposition~\ref{p01} and Lemma~\ref{l12} 
 in the appendix, when the process $u$ satisfies the cyclic 
 condition 
\begin{equation} 
\label{cyclic5} 
 D_{t_1} u_{t_2} ( \omega ) \cdots D_{t_k} u_{t_1} ( \omega ) = 0, 
 \qquad 
 \omega \in \Omega^X, 
 \quad 
 t_1,\ldots ,t_k \in X 
, 
\end{equation} 
 $k \geq 2$, 
 Relation~\eqref{xt} becomes 
\begin{eqnarray} 
\label{xt2} 
\lefteqn{ 
 E_\sigma [ \delta_\sigma ( u )^{n+1} ] 
 = 
 \sum_{k=0}^{n-1} 
 {n \choose k} 
 E_\sigma \left[ 
 \delta_\sigma ( u )^k 
 \int_X 
 u^{n-k+1}_t
 \sigma ( dt ) 
 \right] 
} 
\\ 
\nonumber 
 & & 
 \! \! \! \! \! \! \! \! \! \! \! \! 
 + 
 \sum_{a=0}^n 
 \sum_{j=0}^{a\wedge (b-1)} 
 \sum_{b = a}^n 
 \sum_{{l_0 + \cdots + l_a = n - b 
 \atop l_0 ,\ldots ,l_a \geq 0 
 } 
 } 
 \! \! \! 
 {a \choose j} 
 C_{\eufrak{L}_a,b}^{l_0,n} 
 E_\sigma 
 \left[ 
 \int_{X^{b+1}} 
 \Delta_{\eufrak{s}_{j} } 
 \left( 
 \prod_{q=a+1}^b 
 u_{s_q} 
 \prod_{p=0}^a 
 u^{1+l_p}_{s_p} 
 \right) 
 d\sigma^{b+1} ( \eufrak{s}_b ) 
 \right] 
, 
\end{eqnarray} 
 i.e. the last two terms of \eqref{qas3} vanish when $n=2$. 
 In case $X=\real_+ \times Z$, Condition~\eqref{cyclic5} 
 is satisfied when $u$ is predictable, by the same argument 
 as the one leading to \eqref{cyclic1.111}. 
\\ 
 
 The next Proposition follows from Proposition~\ref{p01} 
 and is used to prove Proposition~\ref{th}. 
\begin{prop} 
\label{c1} 
 Let $N\geq 0$ and let $u\in \lee_{2,1}$ 
 be a bounded process such that 
 $ 
 \displaystyle 
 u \in \bigcap_{p=1}^{N+1} 
 L^\infty ( \Omega^X , L^p_\sigma ( X ) )$ 
 and 
 the integral 
 $\displaystyle \int_X u^n_t \sigma (dt)$ is deterministic, 
 for all $n = 1, \ldots, N+1$, 
 and 
$$ 
 E_\sigma 
 \left[ 
 \int_{X^{a+1}} 
 \Big| 
 \Delta_{s_0} 
 \cdots  
 \Delta_{s_a} 
 \left( 
 \prod_{p=0}^a 
 u^{l_p}_{s_p} 
 \right) 
 \Big| 
 d\sigma^{a+1} ( \eufrak{s}_a ) 
 \right] 
 < \infty, 
$$ 
 $l_0 + \cdots + l_a \leq N+1$, $l_0 ,\ldots ,l_a \geq 1$, 
 $0 \leq a \leq N+1$. 
 Then for all $n = 0,\ldots ,N$ we have 
 $\delta_\sigma ( u ) \in L^{n+1} ( \Omega^X , \pi_\sigma )$ 
 and 
\begin{eqnarray*} 
\lefteqn{ 
 E_\sigma [ \delta_\sigma ( u )^{n+1}] 
 = 
 \sum_{k=0}^{n-1} 
 {n \choose k} 
 \int_X 
 u^{n-k+1}_t
 \sigma ( dt ) 
 E_\sigma \left[ 
 \delta_\sigma ( u )^k 
 \right] 
} 
\\ 
 & & 
 \! \! \! \! \! \! \! \! \! 
 + 
 \! \! \! \! \! 
 \sum_{0\leq j \leq a \leq b \leq n} 
 \sum_{{l_0 + \cdots + l_a = n - b 
 \atop l_0 ,\ldots ,l_a \geq 0 
 } 
 \atop l_{a+1},\ldots ,l_b = 0
 } 
 \! \! \! 
 {a \choose j} 
 C_{\eufrak{L}_a,b}^{l_0,n} 
 E_\sigma 
 \left[ 
 \int_{X^{j+1}} 
 \Delta_{\eufrak{s}_{j} } 
 \prod_{p=0}^j 
 u^{1+l_p}_{s_p} 
 \hskip0.05cm d\sigma^{j+1} ( \eufrak{s}_j ) 
 \right] 
 \prod_{q=j+1}^b 
 \int_X 
 u^{1+l_q}_t 
 \sigma (dt) 
. 
\end{eqnarray*} 
\end{prop} 
\begin{Proof} 
 We apply Proposition~\ref{p01} after integrating 
 in $s_{j+1},\ldots , s_a$ and using \eqref{defdlta}. 
\end{Proof} 
 Consequently if 
 $u : \Omega^X \to \real$ satisfy the hypotheses of Proposition~\ref{c1} 
 and is such that 
\begin{equation} 
\label{lc00000} 
 \int_{X^{j+1}} 
 \Delta_{s_0} 
 \cdots  
 \Delta_{s_j} 
 \left( 
 \prod_{p=0}^j 
 u^{l_p}_{s_p} 
 \right) 
 d\sigma^{j+1} ( \eufrak{s}_j ) 
 = 
 0 
, 
\end{equation} 
 $\pi_\sigma$-a.s., 
 for all $l_0+\cdots + l_j \leq N+1$, 
 $l_0 \geq 1, \ldots , l_j \geq 1$, $j= 1, \ldots ,N$, 
 or simply the cyclic condition 
$$ 
 D_{t_0} u_{t_1} ( \omega ) \cdots D_{t_j} u_{t_0} ( \omega ) = 0, 
 \qquad 
 \omega \in \Omega^X, 
 \quad 
 t_1,\ldots ,t_j \in X 
, 
$$ 
 $j= 1, \ldots ,N$, cf. Lemma~\ref{l12} below, then we have 
$$ 
 E_\sigma [ \delta_\sigma ( u )^{n+1}] 
 = 
 \sum_{k=0}^{n-1} 
 {n \choose k} 
 \int_X 
 u^{n-k+1}_t
 \sigma ( dt ) 
 E_\sigma \left[ 
 \delta_\sigma ( u )^k 
 \right] 
, 
 \qquad 
 n = 0 , \ldots , N, 
$$ 
 i.e. the moments of $\delta_\sigma ( u )$ satisfy the same 
 recurrence relation \eqref{den} 
 as the moments of compensated Poisson integrals. 
\\ 
 
 The next proposition is used to prove Proposition~\ref{p01} 
 with the help of Lemma~\ref{l5}. 
\begin{prop} 
\label{pr1} 
 Let $u:\Omega^X \times X \to \real$ and $v:\Omega^X \times X \to \real$ 
 be bounded processes with compact support in $X$. 
 For all $k \geq 0$ we have 
\begin{eqnarray*} 
\lefteqn{ 
 E_\sigma \left[ 
 \int_X 
 v_s 
 \delta_\sigma ( \varepsilon^+_s u )^k 
 \sigma (ds) 
 \right] 
 = 
 E_\sigma \left[ 
 \delta_\sigma ( u )^k 
 \int_X 
 v_s 
 \sigma (ds) 
 \right] 
} 
\\ 
 & & 
 \! \! \! \! \! \! \! \! \! \! \! \! 
 + 
 \sum_{a=0}^k 
 \sum_{j=0}^a 
 {a \choose j} 
 \sum_{b = a}^k 
 (-1)^{b-a} 
 \! \! \! \! \! \! 
 \sum_{{l_1 + \cdots + l_a = k - b 
 \atop l_1 , \ldots , l_a \geq 0 
 } 
 } 
 C_{\eufrak{L}_a , b } 
 E_\sigma 
 \left[ 
 \int_{X^{b+1}} 
 \Delta_{\eufrak{s}_{j} } 
 \left( 
 v_{s_0} 
 \prod_{q=a+1}^b 
 u_{s_q} 
 \prod_{p=1}^a 
 u^{1+l_p}_{s_p} 
 \right) 
 d\sigma^{b+1} ( \eufrak{s}_b ) 
 \right] 
, 
\end{eqnarray*} 
 where $C_{\eufrak{L}_a , b}$ is defined in \eqref{defin}. 
\end{prop} 
\begin{Proof} 
 This proof is an application of Theorem~\ref{l22} with $F=v_s$. 
 Using Proposition~\ref{prel} below and the expansion 
$$ 
 \prod_{i=0}^j 
 ( I + \Delta_{s_i} ) 
 = 
 \sum_{l=0}^j 
 \sum_{0 \leq i_0 < \cdots < i_l \leq j} 
 \Delta_{s_{i_0}} \cdots \Delta_{s_{i_l}} 
, 
$$ 
 we have, up to the symmetrization 
 due to the integral in $\sigma ( ds_0 ) \cdots \sigma ( ds_a )$ and 
 the summation on $l_1,\ldots ,l_a$, 
\begin{eqnarray*} 
\lefteqn{ 
 \varepsilon^+_{ \eufrak{s}_a } 
 v_{s_0} 
 ( I + D_{s_0} ) 
 \left( 
 \prod_{q=a+1}^b 
 \varepsilon^+_{ \eufrak{s}_a } 
 u_{s_q} 
 \prod_{p=1}^a 
 \varepsilon^+_{ \eufrak{s}_a \setminus s_p } 
 u^{1+l_p}_{s_p} 
 \right) 
} 
\\ 
 & = & 
 \left( 
 \prod_{i=1}^a 
 ( I + \Delta_{s_i} ) 
 \right) 
 \left( 
 v_{s_0} 
 ( I + D_{s_0} ) 
 \left( 
 \prod_{p=1}^a 
 u^{1+l_p}_{s_p} 
 \prod_{q=a+1}^b 
 u_{s_q} 
 \right) 
 \right) 
\\ 
 & = & 
 \left( 
 \prod_{i=1}^a 
 ( I + \Delta_{s_i} ) 
 \right) 
 \left( 
 v_{s_0} 
 \prod_{p=1}^a 
 u^{1+l_p}_{s_p} 
 \prod_{q=a+1}^b 
 u_{s_q} 
 \right) 
\\ 
 & & 
 + 
 \left( 
 \prod_{i=1}^a 
 ( I + \Delta_{s_i} ) 
 \right) 
 \left( 
 v_{s_0} 
 D_{s_0} 
 \left( 
 \prod_{p=1}^a 
 u^{1+l_p}_{s_p} 
 \prod_{q=a+1}^b 
 u_{s_q} 
 \right) 
 \right) 
\\ 
 & = & 
 \left( 
 \prod_{i=1}^a 
 ( I + \Delta_{s_i} ) 
 \right) 
 \left( 
 v_{s_0} 
 \prod_{p=1}^a 
 u^{1+l_p}_{s_p} 
 \prod_{q=a+1}^b 
 u_{s_q} 
 \right) 
\\ 
 & & 
 + 
 \sum_{j=0}^a 
 {a \choose j} 
 \Delta_{s_1} 
 \cdots 
 \Delta_{s_j} 
 \left( 
 v_{s_0} 
 D_{s_0} 
 \left( 
 \prod_{p=1}^a 
 u^{1+l_p}_{s_p} 
 \prod_{q=a+1}^b 
 u_{s_q} 
 \right) 
 \right) 
\\ 
 & = & 
 \varepsilon^+_{ \eufrak{s}_a } 
 v_{s_0} 
 \prod_{q=a+1}^b 
 \varepsilon^+_{ \eufrak{s}_a } 
 u_{s_q} 
 \prod_{p=1}^a 
 \varepsilon^+_{ \eufrak{s}_a \setminus s_p } 
 u^{1+l_p}_{s_p} 
 + 
 \sum_{j=0}^a 
 {a \choose j} 
 \Delta_{s_0} 
 \cdots 
 \Delta_{s_j} 
 \left( 
 v_{s_0} 
 \prod_{p=1}^a 
 u^{1+l_p}_{s_p} 
 \prod_{q=a+1}^b 
 u_{s_q} 
 \right) 
, 
\end{eqnarray*} 
 hence by Theorem~\ref{l22} applied to $G=v_s$ with fixed $s\in X$ 
 we have 
\begin{eqnarray} 
\nonumber 
\lefteqn{ 
 E_\sigma \left[ 
 \int_X v_s 
 \delta_\sigma ( ( I + D_s ) u )^k 
 \sigma ( ds ) 
 \right] 
} 
\\ 
\nonumber 
 & = & 
 \sum_{a=0}^k 
 \sum_{b = 0}^{k-a} 
 (-1)^{b-a} 
 \sum_{l_1 + \cdots + l_a = k - b 
 \atop l_1,\ldots ,l_a \geq 0} 
 C_{\eufrak{L}_a , b} 
\\ 
\nonumber 
 & & 
 E_\sigma \left[  \int_{X^{b+1}} 
 \varepsilon^+_{ \eufrak{s}_a } 
 v_{s_0} 
\right. 
 \left. 
 ( I + D_{s_0} ) 
 \left( 
 \prod_{q=a+1}^b 
 \varepsilon^+_{ \eufrak{s}_a } 
 u_{s_q} 
 \prod_{p=1}^a 
 \varepsilon^+_{ \eufrak{s}_a \setminus s_p } 
 u_{s_p}^{1+l_p} 
 \right) 
 d\sigma^{b+1} ( \eufrak{s}_b ) 
 \right] 
\\ 
\label{1as} 
 & & 
 \! \! \! \! \! \! \! \! \! 
 = 
 \sum_{a=0}^k 
 \sum_{b = 0}^{k-a} 
 (-1)^{b-a} 
 \! \! \! \! \! \sum_{l_1 + \cdots + l_a = k - b 
 \atop l_1,\ldots ,l_a \geq 0} 
 \! \! \! 
 C_{\eufrak{L}_a , b} 
 E_\sigma \left[  \int_{X^{b+1}} 
 \varepsilon^+_{ \eufrak{s}_a } 
 v_{s_0} 
 \prod_{q=a+1}^b 
 \varepsilon^+_{ \eufrak{s}_a } 
 u_{s_q} 
 \prod_{p=1}^a 
 \varepsilon^+_{ \eufrak{s}_a \setminus s_p } 
 u_{s_p}^{1+l_p} 
 \hskip0.05cm d\sigma^{b+1} ( \eufrak{s}_b ) 
 \right] 
\\ 
\nonumber 
 & & 
 + 
 \sum_{a=1}^k 
 \sum_{j=0}^a 
 {a \choose j} 
 \sum_{b = 0}^{k-a} 
 (-1)^{b-a} 
 \sum_{{l_1 + \cdots + l_a = k - b 
 \atop l_1,\ldots ,l_a \geq 0 
 } 
 } 
 C_{\eufrak{L}_a , b} 
\\ 
\nonumber 
 & & 
 E_\sigma 
 \left[ 
 \int_{X^{b+1}} 
 \Delta_{s_0} 
 \cdots 
 \Delta_{s_j} 
 \left( 
 v_{s_0} 
 \prod_{q=a+1}^b 
 u_{s_q} 
 \prod_{p=1}^a 
 u^{1+l_p}_{s_p} 
 \right) 
 d\sigma^{b+1} ( \eufrak{s}_b ) 
 \right] 
\\ 
\nonumber 
 & = & 
 E_\sigma \left[ 
 \delta_\sigma ( u )^k 
 \int_X v_s 
 \hskip0.05cm \sigma ( ds ) 
 \right] 
\\ 
\nonumber 
 & & 
 \! \! \! \! \! \! \! \! \! \! \! 
 + 
 \sum_{a=1}^k 
 \sum_{j=0}^a 
 {a \choose j} 
 \sum_{b = 0}^{k-a} 
 (-1)^{b-a} 
 \! \! \! \! \! 
 \sum_{{l_1 + \cdots + l_a = k - b 
 \atop l_1,\ldots ,l_a \geq 0 
 }
 } 
 \! \!\! \! 
 C_{\eufrak{L}_a , b} 
 E_\sigma 
 \left[ 
 \int_{X^{b+1}} 
 \Delta_{\eufrak{s}_{j} } 
 \left( 
 v_{s_0} 
 \prod_{q=a+1}^b 
 u_{s_q} 
 \prod_{p=1}^a 
 u^{1+l_p}_{s_p} 
 \right) 
 d\sigma^{b+1} ( \eufrak{s}_b ) 
 \right] 
, 
\end{eqnarray} 
 where we identified 
 $\displaystyle 
 E_\sigma \left[ 
 \delta_\sigma ( u )^k 
 \int_X v_s 
 \sigma (ds) 
 \right] 
 $ 
 to \eqref{1as} 
 on the last step, 
 by another application of Theorem~\ref{l22} to $F=\int_X v_s \sigma (ds)$. 
\end{Proof} 
\section{Appendix} 
\label{app2} 
 In this appendix we state some combinatorial results that have 
 been used above. 
\begin{prop} 
\label{prel} 
 Let $u : \Omega^X \times X \to \real$ be a measurable process. 
 For all $0 \leq j , p \leq n$ we have the relation 
\begin{equation} 
\label{therel} 
 \prod_{p=0}^n 
 \varepsilon^+_{ \eufrak{s}_j \setminus s_p } 
 u_{s_p} 
 = 
 \left( 
 \prod_{i=0}^j 
 ( I + \Delta_{s_i} ) 
 \right) 
 \prod_{p=0}^n 
 u_{s_p} 
, 
\end{equation} 
 for mutually different $\eufrak{s}_n = ( s_0 , \ldots , s_n ) \subset X$. 
\end{prop} 
\begin{Proof} 
 We will prove Relation~\eqref{therel} for all $n\geq 0$ 
 by induction on $j \in \{ 0 , 1, \ldots , n \}$. 
 Clearly for $j=0$ the relation holds since 
\begin{eqnarray*} 
 u_{s_0} 
 \prod_{p=1}^n 
 \varepsilon^+_{s_0} 
 u_{s_p} 
 & = & 
 u_{s_0} 
 \prod_{p=1}^n 
 u_{s_p} 
 + 
 u_{s_0} 
 \sum_{ 
 { 
 \Xi_1 \cup \cdots \cup \Xi_n = \{ s_0 \} 
 } 
 } 
 D_{\Xi_1} 
 u_{s_1} 
 \cdots 
 D_{\Xi_n} 
 u_{s_n} 
\\ 
 & = & 
 u_{s_0} 
 \prod_{p=1}^n 
 u_{s_p} 
 + 
 \sum_{ 
 { 
 \Xi_0 \cup \cdots \cup \Xi_n = \{ s_0 \} 
 \atop 
 s_0 \notin \Xi_0 
 } 
 } 
 D_{\Xi_0} 
 u_{s_0} 
 \cdots 
 D_{\Xi_n} 
 u_{s_n} 
\\ 
 & = & 
 ( I + \Delta_{s_0} ) 
 \prod_{p=0}^n 
 u_{s_p} 
. 
\end{eqnarray*} 
 Next, assuming that \eqref{therel} holds at the 
 rank $j \in \{ 0,1,\ldots ,n-1 \}$ and taking 
 $\{ s_0 , \ldots , s_n \} \subset X$ mutually different 
 we have 
\begin{eqnarray*} 
\lefteqn{ 
 \prod_{p=0}^n 
 \varepsilon^+_{ \eufrak{s}_{j+1} \setminus s_p } 
 u_{s_p} 
 =
 \prod_{p=0}^n 
 \left( 
 ( I + 
 {\bf 1}_{\{ p \not= j+1 \} } 
 D_{s_{j+1}} ) 
 \prod_{i=0 \atop i\not= p}^j 
 ( I + D_{s_i} ) 
 u_{s_p} 
 \right) 
} 
\\ 
 & = & 
 \sum_{ 
 { 
 \Xi_0 \cup \cdots \cup \Xi_n \subset \{ s_{j+1} \} 
 \atop 
 s_{j+1} \notin \Xi_{j+1} 
 } 
 } 
 \prod_{p=0}^n 
 \left( 
 \prod_{i=0 \atop i\not= p}^j 
 ( I + D_{s_i} ) 
 D_{\Xi_p} 
 u_{s_p} 
 \right) 
\\ 
 & = & 
 \sum_{ 
 { 
 \Xi_0 \cup \cdots \cup \Xi_n \subset \{ s_{j+1} \} 
 \atop 
 s_{j+1} \notin \Xi_{j+1} 
 } 
 } 
 \left( 
 \prod_{i=0}^j 
 ( I + \Delta_{s_i} ) 
 \right) 
 \prod_{p=0}^n 
 D_{\Xi_p} 
 u_{s_p} 
\\ 
 & = & 
 \sum_{ 
 { 
 \Xi_0 \cup \cdots \cup \Xi_n \subset \{ s_{j+1} \} 
 \atop 
 s_{j+1} \notin \Xi_{j+1} 
 } 
 } 
 \sum_{ 
 { 
 \Theta_0 \cup \cdots \cup \Theta_n = \{ s_0 , s_1 , \ldots , s_j \} 
 \atop 
 s_0 \notin \Theta_0 , \ldots , s_j \notin \Theta_j 
 } 
 } 
 D_{\Theta_0} 
 D_{\Xi_0} 
 u_{s_0} 
 \cdots 
 D_{\Theta_n} 
 D_{\Xi_n} 
 u_{s_n} 
\\ 
 & = & 
 \sum_{ 
 { 
 \Theta_0 \cup \cdots \cup \Theta_n = \{ s_0 , s_1 , \ldots , s_j \} 
 \atop 
 s_0 \notin \Theta_0 , \ldots , s_j \notin \Theta_j 
 } 
 } 
 D_{\Theta_0} 
 u_{s_0} 
 \cdots 
 D_{\Theta_n} 
 u_{s_n} 
\\ 
 & & 
 + 
 \sum_{ 
 { 
 \Xi_0 \cup \cdots \cup \Xi_n = \{ s_{j+1} \} 
 \atop 
 s_{j+1} \notin \Xi_{j+1} 
 } 
 } 
 \sum_{ 
 { 
 \Theta_0 \cup \cdots \cup \Theta_n = \{ s_0 , s_1 , \ldots , s_j \} 
 \atop 
 s_0 \notin \Theta_0 , \ldots , s_j \notin \Theta_j 
 } 
 } 
 D_{\Xi_0} 
 D_{\Theta_0} 
 u_{s_0} 
 \cdots 
 D_{\Xi_n} 
 D_{\Theta_n} 
 u_{s_n} 
\\ 
 & = & 
 \left( 
 \prod_{i=0}^j 
 ( I + \Delta_{s_i} ) 
 \right) 
 \prod_{p=0}^n 
 u_{s_p} 
 + 
 \Delta_{s_{j+1}} 
 \left( 
 \prod_{i=0}^j 
 ( I + \Delta_{s_i} ) 
 \right) 
 \prod_{p=0}^n 
 u_{s_p} 
\\ 
 & = & 
 \left( 
 \prod_{i=0}^{j+1} 
 ( I + \Delta_{s_i} ) 
 \right) 
 \prod_{p=0}^n 
 u_{s_p} 
. 
\end{eqnarray*} 
\end{Proof} 
 Finally in the next lemma, which is used to prove 
 Corollary~\ref{c09}, we show that Relation~\eqref{lc00000} 
 is satisfied provided $D_su_t ( \omega )$ satisfies the cyclic condition 
 \eqref{cyclic2}. 
\begin{lemma} 
\label{l12} 
 Let $N\geq 1$ and assume that $u : \Omega^X \times X \to \real$ 
 satisfies the cyclic condition 
\begin{equation} 
\label{cyclic2} 
 D_{t_0} u_{t_1} ( \omega ) \cdots D_{t_j} u_{t_0} ( \omega ) = 0, 
 \qquad 
 \omega \in \Omega^X, 
 \quad 
 t_0, t_1, \ldots ,t_j \in X 
, 
\end{equation} 
 for $j = 1, \ldots , N$. 
 Then we have 
$$ 
 \Delta_{t_0} 
 \cdots  
 \Delta_{t_j} 
 \left( 
 u_{t_0} ( \omega ) \cdots u_{t_j} ( \omega ) 
 \right) 
 = 0, 
 \qquad 
 \omega \in \Omega^X, 
 \quad 
 t_0, t_1, \ldots ,t_j \in X 
, 
$$ 
 for $j = 1, \ldots , N$. 
\end{lemma} 
\begin{Proof} 
 By Definition~\ref{d1} we have 
\begin{equation} 
\label{fjkl}
 \Delta_{t_0} 
 \cdots  
 \Delta_{t_j} 
 \prod_{p=0}^j 
 u_{t_p} 
 = 
 \sum_{ 
 { 
 \Theta_0 \cup \cdots \cup \Theta_j = \{ t_0 , t_1 , \ldots , t_j \} 
 \atop 
 t_0 \notin \Theta_0 , \ldots , t_j \notin \Theta_j 
 } 
 } 
 D_{\Theta_0} 
 u_{t_0} 
 \cdots 
 D_{\Theta_j} 
 u_{t_j} 
, 
\end{equation}  
 $t_0,\ldots ,t_j \in X$, $j = 2 , \ldots , N$. 
 Without loss of generality we may assume that 
 $\{ t_0, t_1,\ldots ,t_j \}$ are not equal to eachother and 
 that $\Theta_0 \not= \emptyset, \ldots , \Theta_j \not= \emptyset$ 
 and $\Theta_k \cap \Theta_l = \emptyset$, $0 \leq k \not=l \leq j$ 
 in the above sum. 
 In this case we can construct a sequence $(k_1,\ldots ,k_i)$ by 
 choosing 
$$ 
 t_0 \not= t_{k_1} \in \Theta_0, \ t_{k_2} \in \Theta_{k_1}, \ 
 \ldots, 
 t_{k_{i-1}} \in \Theta_{k_{i-2}}, 
$$ 
 until $t_{k_i} = t_0 \in \Theta_{k_{i-1}}$ for some $i\in \{2,\ldots ,j \}$ 
 since $\Theta_0 \cap \cdots \cap \Theta_j =\emptyset$ 
 and 
 $\Theta_0 \cup \cdots \cup \Theta_j = \{ t_0 , t_1 , \ldots , t_j \}$. 
 Hence by \eqref{cyclic2} we have 
$$ 
 D_{t_{k_1}} 
 u_{t_0} 
 D_{t_{k_2}} 
 u_{t_{k_1}} 
 \cdots 
 D_{t_{k_{i-1}}} 
 u_{t_{k_{i-2}}} 
 D_{t_{k_0}} 
 u_{t_{k_{i-1}}} 
 = 0 
, 
$$ 
 by \eqref{cyclic2}, which implies 
$$ 
 D_{\Theta_0} 
 u_{t_0} 
 D_{\Theta_{k_1}} 
 u_{t_{k_1}} 
 \cdots 
 D_{\Theta_{k_{i-1}}} 
 u_{t_{k_{i-1}}} 
 = 0 
, 
$$ 
 since 
$$ 
 ( t_{k_1} , t_{k_2} , \ldots , t_{k_{i-1}} , t_0 ) \in 
 \Theta_0 \times \Theta_{k_1} \times \cdots \times \Theta_{k_{i-1}} 
, 
$$ 
 hence 
$$ 
 D_{\Theta_0} 
 u_{t_0} 
 D_{\Theta_{k_1}} 
 u_{t_{k_1}} 
 \cdots 
 D_{\Theta_{k_j}} 
 u_{t_j} 
 = 0 
, 
$$ 
 and \eqref{fjkl} vanishes. 
\end{Proof} 
 Again, in case $X=\real_+$, Condition~\eqref{cyclic2} holds 
 in particular when either 
$ 
 D_s u_t 
 = 
 0, 
$ 
 $0 \leq s \leq t$, 
 as in \eqref{erf1.111}, 
 resp.  
$ 
 D_t u_s 
 = 
 0$, 
 $0 \leq s \leq t$, which is the case when $u$ is backward, 
 resp. forward, predictable. 
 
\footnotesize 

\def\cprime{$'$} \def\polhk#1{\setbox0=\hbox{#1}{\ooalign{\hidewidth
  \lower1.5ex\hbox{`}\hidewidth\crcr\unhbox0}}}
  \def\polhk#1{\setbox0=\hbox{#1}{\ooalign{\hidewidth
  \lower1.5ex\hbox{`}\hidewidth\crcr\unhbox0}}} \def\cprime{$'$}

\end{document}